%% file: BCD-G-SIMAX.tex
\documentclass[final,onefignum,onetabnum]{siamart190516}

\usepackage{braket,amsfonts}

\usepackage{bigdelim}

\usepackage{array}

\usepackage{pgfplots}

\newsiamthm{claim}{Claim}
\newsiamremark{remark}{Remark}
\newsiamremark{hypothesis}{Hypothesis}
\Crefname{hypothesis}{Hypothesis}{Hypotheses}

\usepackage{algorithmic}
\usepackage{graphicx,epstopdf}

\Crefname{ALC@unique}{Line}{Lines}

\usepackage{amsopn}

\usepackage{xspace}
\usepackage{bold-extra}
\usepackage[most]{tcolorbox}

\colorlet{texcscolor}{blue!50!black}
\colorlet{texemcolor}{red!70!black}
\colorlet{texpreamble}{red!70!black}
\colorlet{codebackground}{black!25!white!25}

\usepackage{amsfonts}

\newcommand{\overbar}[1]{\mkern 1.5mu\overline{\mkern-1.5mu#1\mkern-1.5mu}\mkern 1.5mu}

\lstdefinestyle{siamlatex}{%
style=tcblatex,
texcsstyle=*\color{texcscolor},
texcsstyle=[2]\color{texemcolor},
keywordstyle=[2]\color{texemcolor},
moretexcs={cref,Cref,maketitle,mathcal,text,headers,email,url},
}

\tcbset{%
colframe=black!75!white!75,
coltitle=white,
colback=codebackground, % bottom/left side
colbacklower=white, % top/right side
fonttitle=\bfseries,
arc=0pt,outer arc=0pt,
top=1pt,bottom=1pt,left=1mm,right=1mm,middle=1mm,boxsep=1mm,
leftrule=0.3mm,rightrule=0.3mm,toprule=0.3mm,bottomrule=0.3mm,
listing options={style=siamlatex}
}

\newtcblisting[use counter=example]{example}[2][]{%
title={Example~\thetcbcounter: #2},#1}

\newtcbinputlisting[use counter=example]{\examplefile}[3][]{%
title={Example~\thetcbcounter: #2},listing file={#3},#1}

\DeclareTotalTCBox{\code}{ v O{} }
{ %fontupper=\ttfamily\color{texemcolor},
fontupper=\ttfamily\color{black},
nobeforeafter,
tcbox raise base,
colback=codebackground,colframe=white,
top=0pt,bottom=0pt,left=0mm,right=0mm,
leftrule=0pt,rightrule=0pt,toprule=0mm,bottomrule=0mm,
boxsep=0.5mm,
#2}{#1}

\patchcmd\newpage{\vfil}{}{}{}
\usepackage{amssymb}
\definecolor{cobalt}{rgb}{0.0, 0.28, 0.67}
\definecolor{darkblue}{rgb}{0.0, 0.0, 0.65}

\definecolor{darkred}{rgb}{0.8,0,0}
\definecolor{darkgreen}{rgb}{0,0.46,0}
\definecolor{purple}{rgb}{0.6,0,0.5}

\usepackage{tikz}
\usetikzlibrary{shapes.geometric, arrows}

\tikzstyle{startstop} = [rectangle, rounded corners, minimum width=2.4cm, minimum height=0.65cm,text centered, draw=black, fill=green!30]
\tikzstyle{startstop0} = [rectangle, rounded corners, minimum width=2.4cm, minimum height=0.65cm,text centered, draw=white, fill=white!30]
\tikzstyle{io} = [trapezium, trapezium left angle=70, trapezium right angle=110, minimum width=3cm, minimum height=1cm, text centered, draw=black, fill=blue!30]
\tikzstyle{process} = [rectangle, minimum width=3cm, minimum height=1cm, text centered, draw=black, fill=orange!30]
\tikzstyle{decision} = [diamond, minimum width=3cm, minimum height=1cm, text centered, draw=black, fill=green!30]
\tikzstyle{arrow} = [thick,->,>=stealth]

\newcommand{\matr}[1]{\boldsymbol{#1}}

\newcommand{\vect}[1]{\boldsymbol{#1}}

\newcommand{\T}{{\sf T}} % transposition
\renewcommand{\H}{{\sf H}} % Hermitian transposition

\makeatletter
\newcommand{\rank}[1]{\mathop{\operator@font rank}\{#1\}}
\newcommand{\colrank}[1]{\mathop{\operator@font colrank}\{#1\}}
\newcommand{\krank}[1]{\mathop{\operator@font krank}\{#1\}}
\newcommand{\trace}[1]{\mathop{\operator@font tr}(#1)}
\newcommand{\symmm}[1]{\mathop{\operator@font sym}(#1)}
\newcommand{\skeww}[1]{\mathop{\operator@font skew}(#1)}
\newcommand{\Diag}[1]{\mathop{\operator@font Diag}\{#1\}} % a diagonal matrix
\newcommand{\diag}[1]{\mathop{\operator@font diag}\{#1\}} % a vector
\newcommand{\Span}[1]{\mathop{\operator@font Span}\{#1\}} % a space
\newcommand{\argmin}{\mathop{\operator@font argmin}}

\newcommand{\offdiag}[1]{\mathop{\operator@font offdiag}\{#1\}} % a vector
\newcommand{\Proj}[2]{\mathop{\operator@font Proj_{#1}}{#2}}
\newcommand{\ProjGrad}[2]{\mathop{{\operator@font grad} }#1(#2)}
\makeatother

\newcommand{\eqdef}{\stackrel{\sf def}{=}}
\newcommand{\RR}{\mathbb{R}}
\newcommand{\CC}{\mathbb{C}}
\newcommand{\NN}{\mathbb{N}}

\newcommand{\Gmat}[3]{\matr{Q}^{(#1,#2,#3)}}

\newcommand{\R}{\Re}
\newcommand{\I}{\Im}
\newcommand{\Tang}[2]{\mathbf{T}_{#1}{#2}}

\newcommand{\Utwo}{\Psi}
\newcommand{\hij}[3]{h^{(Q)}_{(#1,#2), #3}}
\newcommand{\Gamij}[3]{\matr{\Gamma}^{(#1,#2,#3)}}
\newcommand{\RInner}[2]{\left\langle #1, #2\right\rangle_{\Re}} % Inner product

\makeatletter
\def\bbordermatrix#1{\begingroup \m@th
\@tempdima 4.75\p@
\setbox\z@\vbox{%
\def\cr{\crcr\noalign{\kern2\p@\global\let\cr\endline}}%
\ialign{$##$\hfil\kern2\p@\kern\@tempdima&\thinspace\hfil$##$\hfil
&&\quad\hfil$##$\hfil\crcr
\omit\strut\hfil\crcr\noalign{\kern-\baselineskip}%
#1\crcr\omit\strut\cr}}%
\setbox\tw@\vbox{\unvcopy\z@\global\setbox\@ne\lastbox}%
\setbox\tw@\hbox{\unhbox\@ne\unskip\global\setbox\@ne\lastbox}%
\setbox\tw@\hbox{$\kern\wd\@ne\kern-\@tempdima\left[\kern-\wd\@ne
\global\setbox\@ne\vbox{\box\@ne\kern2\p@}%
\vcenter{\kern-\ht\@ne\unvbox\z@\kern-\baselineskip}\,\right]$}%
\null\;\vbox{\kern\ht\@ne\box\tw@}\endgroup}
\makeatother

\usepackage{amsmath}
\usepackage{enumerate}
\usepackage{xcolor}
\usepackage[mathscr]{euscript}
\usepackage{url}
\usepackage[all]{xy}
\usepackage{booktabs}
\usepackage{graphicx}
\usepackage{epstopdf}
\usepackage{algorithmic}
\usepackage{mathtools}
\usepackage{makecell}
\usepackage{multirow}
\usepackage[ruled,vlined]{algorithm2e}

\crefname{algocf}{alg.}{algs.}
\Crefname{algocf}{Algorithm}{Algorithms}

\newtheorem{example1}[theorem]{Example}
\newtheorem{setting}[theorem]{Update rule}

\newcommand{\matrelem}[1]{\mathrm{#1}}

\renewcommand{\H}{{\sf H}}      % Hermitian transposition
\renewcommand{\TH}{{\blacklozenge}} 
\newcommand{\expp}[1]{\mathrm{exp}\left(#1\right)}

\newcommand{\UNN}[1]{\mathbf{U}_{#1}}
\newcommand{\SUNN}[1]{\mathbf{SU}_{#1}}

\newcommand{\GL}{\mathbf{GL}}
\newcommand{\SL}{\mathbf{SL}}
\newcommand{\tr}{\mbox{tr}}
\newcommand{\St}{\mathbf{St}}

\newcommand{\UT}{\mathbf{UT}}
\newcommand{\SUT}{\mathbf{SUT}}
\newcommand{\EUT}{\mathbf{EUT}}

\newcommand{\LT}{\mathbf{LT}}
\newcommand{\SLT}{\mathbf{SLT}}
\newcommand{\ELT}{\mathbf{ELT}}

\newcommand{\DT}{\mathbf{D}}

\newcommand{\sll}{\mathfrak{sl}}

\newcommand{\sut}{\mathfrak{sut}}
\newcommand{\eut}{\mathfrak{eut}}

\newcommand{\slt}{\mathfrak{slt}}
\newcommand{\su}{\mathfrak{su}}

\newcommand{\sdt}{\mathfrak{d}}
 
\newcommand{\Lmat}[3]{\matr{L}^{(#1,#2,#3)}}
\newcommand{\Umat}[3]{\matr{U}^{(#1,#2,#3)}}
\newcommand{\Dmat}[3]{\matr{D}^{(#1,#2,#3)}}
 
\newcommand{\ue}{\mathrm{e}}
\newcommand{\ui}{i}

\graphicspath{{./matlab/}}

\begin{tcbverbatimwrite}{tmp_\jobname_header.tex}
\title{Convergence of gradient-based block coordinate descent algorithms for non-orthogonal joint approximate diagonalization of matrices\thanks{Submitted to the editors on Nov. 3, 2021; revised July 1st, 2022; revised Nov. 25, 2022. 
\funding{This work was supported in part by the National Natural Science Foundation of China (No. 11601371), the Guangdong Basic and Applied Basic Research Foundation (No. 2021A1515010232), and Agence Nationale de Recherche (ANR-19-CE23-0021).}}}

\author{Jianze Li\thanks{Shenzhen Research Institute of Big Data, The Chinese University of Hong Kong, Shenzhen, China (\email{lijianze@gmail.com}).}
\and Konstantin Usevich\thanks{Universit\'{e} de Lorraine, CNRS, CRAN, Nancy, France (\email{konstantin.usevich@cnrs.fr}).} \and Pierre Comon\thanks{Univ. Grenoble Alpes, CNRS, Grenoble INP, GIPSA-Lab, France (\email{pierre.comon@gipsa-lab.fr}).}}

\headers{Gradient-based block coordinate descent algorithms}{Jianze Li, Konstantin Usevich and Pierre Comon}
\end{tcbverbatimwrite}
\input{tmp_\jobname_header.tex}

\ifpdf
\hypersetup{ pdftitle={Guide to Using SIAM'S \LaTeX\ Style} }
\fi

\begin{document}
\maketitle

%% ------------------------------------------------------------------
%% ABSTRACT
%% ------------------------------------------------------------------
\begin{tcbverbatimwrite}{tmp_\jobname_abstract.tex}
\begin{abstract}
In this paper, we propose a gradient-based block coordinate descent (BCD-G) framework to solve the joint approximate diagonalization of matrices defined on the product of the complex Stiefel manifold and the special linear group.
Instead of the cyclic fashion, 
we choose a block optimization based on the Riemannian gradient. 
To update the first block variable in the complex Stiefel manifold, we use the well-known line search descent method.
To update the second block variable in the special linear group, based on four kinds of different elementary transformations, we construct three classes: GLU, GQU and GU, and then get three BCD-G algorithms: BCD-GLU, BCD-GQU and BCD-GU. 
We establish the global and weak convergence of these three algorithms using the \L{}ojasiewicz gradient inequality under the assumption that the iterates are bounded. We also propose a gradient-based Jacobi-type framework to solve the joint approximate diagonalization of matrices defined on the special linear group. As in the BCD-G case, using the GLU and GQU classes of elementary transformations, we focus on the Jacobi-GLU and Jacobi-GQU algorithms and establish their global and weak convergence. 
All the algorithms and convergence results described in this paper also apply to the real case. 
\end{abstract}


\end{tcbverbatimwrite}
\input{tmp_\jobname_abstract.tex}
%% ------------------------------------------------------------------
%% END HEADER
%% ------------------------------------------------------------------

\section{Introduction}
Let $1\leq m\leq n$. 
Given a complex matrix $\matr{Z}\in\CC^{n\times m}$, we denote by $\matr{Z}^{\T}$, $\matr{Z}^{*}$ and $\matr{Z}^{\H}$ its \emph{transpose, conjugate} and \emph{conjugate transpose,} respectively. 
We shall also use  $(\cdot)^{\TH}$ to denote either $(\cdot)^{\T}$ or $(\cdot)^{\H}$.
A complex matrix $\matr{A}\in\CC^{n\times n}$ is called \emph{Hermitian} if $\matr{A}^{\H}=\matr{A}$.
It is called \emph{complex symmetric} if $\matr{A}^{\T}=\matr{A}$. 
Let $\{\matr{A}^{(\ell)}\}_{1\leq\ell\leq \matrelem{L}}\subseteq\CC^{n\times n}$ be a set of complex matrices. 
The well-known \emph{blind source separation} (BSS) problem \cite{Como94:sp,Como10:book,maurandi2014decoupled,souloumiac2009nonorthogonal} can be formulated as finding a full column rank matrix $\matr{Z}\in\CC^{n\times m}$ to make the matrices $\matr{W}^{(\ell)}=\matr{Z}^{\TH}\matr{A}^{(\ell)} \matr{Z}\in\CC^{m\times m}$ simultaneously as diagonal as possible. 
A natural idea is to solve the \emph{joint approximate diagonalization of matrices} (JADM) problem, which consists in minimizing
\begin{equation}\label{eq:sym_tensor_diagonalization-gener-0}
f(\matr{Z}) = \sum_{\ell=1}^{\matrelem{L}} \|\offdiag{\matr{W}^{(\ell)}}\|^2,
\end{equation}
where $\matr{Z}\in\CC^{n\times m}$ is a full column rank matrix, and $\offdiag{\cdot}$ is the \emph{zero diagonal} operator, setting all the diagonal elements of a square matrix in $\CC^{m\times m}$ to zero.

Note that, the set of full-column rank matrices is not closed (the limit of a sequence of full column rank matrices can be rank deficient), and therefore problem \cref{eq:sym_tensor_diagonalization-gener-0} is ill-posed.
For example,  for a full column rank matrix $\matr{Z}\in\CC^{n\times m}$ and nonzero $\lambda\in\CC$, we have $\lim_{\lambda\to 0}f(\lambda\matr{Z}) = \lim_{\lambda\to 0} |\lambda|^4f(\matr{Z}) =0$. 
To tackle this issue, it is first necessary to use scale- and permutation-invariant cost functions \cite{afsari2004gradient,yeredor2002non}.
Second, the set of matrices $\matr{Z}$ must be restricted to a smaller \emph{closed} subset. Several possibilities can be envisaged, \emph{e.g.}, a restriction to the special linear group $\SL_{m}(\CC)$ in the square case $m=n$. In this paper, we follow the latter approach as it will be discussed later. 

Problem \cref{eq:sym_tensor_diagonalization-gener-0} has been widely used in BSS  and \emph{Independent component analysis} (ICA)  \cite{chabriel2014joint,Como10:book,Afsa07:honolulu,andre2018new}, and has the following well-known special cases:
\begin{itemize}
\item \emph{joint approximate diagonalization of Hermitian matrices} (JADM-H) \cite{sorensen2009approximate,maurandi2014decoupled}: $(\cdot)^{\TH}=(\cdot)^{\H}$, $\matr{A}^{(\ell)}  \in\CC^{n\times n}$ is Hermitian for $1\leq\ell\leq \matrelem{L}$;
\item \emph{joint approximate diagonalization of complex symmetric matrices} (JADM-CS) \cite{maurandi2014decoupled}: $(\cdot)^{\TH}=(\cdot)^{\T}$, $\matr{A}^{(\ell)} \in\CC^{n\times n}$ is complex symmetric for $1\leq\ell\leq \matrelem{L}$;
\item \emph{joint approximate diagonalization of real symmetric matrices} (JADM-RS) \cite{afsari2004gradient,afsari2006simple}: over real $\matr{Z}$, $(\cdot)^{\TH}=(\cdot)^{\T}$, 
$\matr{A}^{(\ell)} \in\RR^{n\times n}$ is real symmetric for $1\leq\ell\leq \matrelem{L}$. 
\end{itemize}

Many classic approaches use prewhitening to reduce the problem \cref{eq:sym_tensor_diagonalization-gener-0} to orthogonal (and square) diagonalization case \cite{Cardoso93:JADE,cardoso1996jacobi,Como94:sp,LUC2017globally,LUC2018,li2020convergence,li2019jacobi,ULC2019}.
This, however, results in a two-step procedure, which may not be optimal in the statistical sense and may suffer more from noise. 
Therefore, the non-orthogonal joint diagonalization attracted considerable interest in the literature. 
In particular, to solve the JADM-RS problem, 
Jacobi-type algorithms were introduced based on the LU and QR decompositions in \cite{afsari2006simple}, and on the Givens transformations, hyperbolic transformations, and diagonal transformations in \cite[Eq. (9)]{souloumiac2009nonorthogonal}.  
To solve the JADM-H problem, Jacobi-type algorithms were proposed based on the LU decomposition in \cite{maurandi2014decoupled,maurandi2014jacobi}, and based on the QL decomposition in \cite{sorensen2009approximate}.
To solve the JADM-CS problem, a Jacobi-type algorithm was proposed based on the LU decomposition in \cite{maurandi2013fast,maurandi2014decoupled}. 
However, to our knowledge, there was no theoretical result about the convergence of these Jacobi-type algorithms in the literature.
In addition, mostly the square ($m=n$) case was considered. 

In this paper, we consider the general rectangular case of  \cref{eq:sym_tensor_diagonalization-gener-0}, with $\matr{Z}$ restricted to a $\SL_{m}(\CC)$-like subset. 
By using a reformulation of the problem, we develop optimization algorithms on manifolds, and provide convergence results.
An overview of the contributions is provided in the rest of the section.

\subsection{Search space and reformulations of the problem}
Let $\GL_{m}(\CC)\eqdef\{\matr{X}\in\CC^{m\times m}, {\rm det}(\matr{X})\neq 0\}$ (resp.  $\SL_{m}(\CC)\eqdef\{\matr{X}\in\GL_{m}(\CC), {\rm det}(\matr{X})=1\}$) be the \emph{general} (resp. \emph{special}) \emph{linear
group}. 
We define the \emph{rectangular special linear set} as
\begin{align}\label{eq:def_rsl}
\textbf{RSL}(m,n,\CC) \eqdef \{\matr{Z}\in\CC^{n\times m}, \matr{Z}^{\H}\matr{Z}\in\SL_{m}(\CC)\}.
\end{align}
Every matrix in $\textbf{RSL}(m,n,\CC)$ is of full column rank, and, moreover this set is closed. 
Thus the problem of rank deficiency or trivial solution at $\matr{0}$ does not appear when optimizing \cref{eq:sym_tensor_diagonalization-gener-0} over $\textbf{RSL}(m,n,\CC)$, since $\lambda\matr{Z}\notin\textbf{RSL}(m,n,\CC)$ if $|\lambda|\neq1$.
Still, this remains a difficult optimization problem, since the feasible region $\textbf{RSL}(m,n,\CC)$ is neither convex nor compact, and the function $f(\matr{Z})$ is a  quartic polynomial. 
In what follows, we provide a reformulation of the problem for two scenarios.

\begin{itemize}
\item \emph{General (rectangular) case.}
Let 
$\St(m,n,\CC) \eqdef\{\matr{Y}\in\CC^{n\times m}, \matr{Y}^{\H}\matr{Y}= \matr{I}_m\}$ 
be the \emph{complex Stiefel manifold}. 
We have the following simple result: 
\begin{lemma}\label{lem:rela-equiavalent}
A complex matrix $\matr{Z}\in\textbf{RSL}(m,n,\CC)$ if and only if there exist $\matr{Y}\in\St(m,n,\CC)$ and $\matr{X}\in\SL_{m}(\CC)$ such that $\matr{Z}=\matr{Y}\matr{X}$.
\end{lemma} 
By \Cref{lem:rela-equiavalent}, problem \cref{eq:sym_tensor_diagonalization-gener-0} over $\textbf{RSL}(m,n,\CC)$ is equivalent to minimizing 
\begin{equation}\label{eq:sym_tensor_diagonalization-gener-0-equi}
f: \St(m,n,\CC)\times\SL_{m}(\CC)\rightarrow\RR^{+},\ \ (\matr{Y},\matr{X})\mapsto \sum_{\ell=1}^{\matrelem{L}} \|\offdiag{\matr{W}^{(\ell)}}\|^2,
\end{equation}
where $\matr{W}^{(\ell)}=(\matr{Y}\matr{X})^{\TH}\matr{A}^{(\ell)}(\matr{Y}\matr{X})\in\CC^{m\times m}$.

\item  \emph{Square case (second reformulation).}
This is a special case of 
\cref{eq:sym_tensor_diagonalization-gener-0-equi}, when we assume $\matr{Y}_{*}\in\St(m,n,\CC)$  to be fixed (for example, it is found in advance by some other method, \emph{e.g.}, PCA  \cite{Como92:elsevier,Como94:sp,Como10:book}, which is a common procedure for dimensionality and noise reduction). 
Denote  $\matr{B}^{(\ell)}=\matr{Y}_{*}^{\TH}\matr{A}^{(\ell)} \matr{Y}_{*}$ for $1\leq\ell\leq \matrelem{L}$.  
Then the cost function  \cref{eq:sym_tensor_diagonalization-gener-0-equi} becomes  
\begin{equation}\label{eq:sym_tensor_diagonalization-gener-1}
g: \SL_{m}(\CC)\rightarrow\RR^{+},\ \ \matr{X}\mapsto  g(\matr{X})=  \sum_{\ell=1}^{\matrelem{L}} \|\offdiag{\matr{W}^{(\ell)}}\|^2,
\end{equation}
where $\matr{W}^{(\ell)}=\matr{X}^{\TH}\matr{B}^{(\ell)} \matr{X}\in\CC^{m\times m}$. 
Alternatively, this case may appear when $m=n$ in \cref{eq:sym_tensor_diagonalization-gener-0}.
Indeed,  $\textbf{RSL}(m,m,\CC) = \{\matr{Z}\in\CC^{m\times m}, |\det(\matr{Z})| = 1\}$,
and since \cref{eq:sym_tensor_diagonalization-gener-0} is invariant with respect to multiplication by a unimodular scalar, we can optimize it over  $\SL_{m}(\CC)$ instead. 
\end{itemize}

\subsection{Contributions}

In this paper, 
to solve problem \cref{eq:sym_tensor_diagonalization-gener-0-equi}, which is defined on the product of $\St(m,n,\CC)$ and $\SL_{m}(\CC)$, the \emph{gradient-based block coordinate descent} (BCD-G) algorithms (\Cref{alg:BCD-G-RL}) will be proposed in \Cref{subsecc:bcd-g} (more detailedly in \Cref{subsec:Jaco_G}), which chooses a block  optimization based on the Riemannian gradient.
This is similar to the gradient-based way of choosing index pairs in the Jacobi-G algorithms on the orthogonal group \cite{IshtAV13:simax,LUC2017globally} or unitary group \cite{ULC2019}. 
Then their \emph{global convergence}\footnote{For any starting point, the iterates converge to a limit point as a whole sequence.} and \emph{weak convergence}\footnote{Every accumulation point is a stationary point, \emph{i.e.}, the Riemannian gradient is equal to 0.} will be established in \Cref{sec:BCD-G} using the \emph{\L{}ojasiewicz gradient inequality} \cite{loja1965ensembles,lojasiewicz1993geometrie,absil2005convergence,Usch15:pjo}, under the assumption that the iterates $\matr{\omega}_{k}$ are bounded, that is, 
there exists a universal positive constant  $\mathrm{M}_{\omega}>0$ such that 
\begin{equation}\label{eq:condi_bounded_X_gene}
\|\matr{\omega}_{k}\|\leq \mathrm{M}_{\omega}
\end{equation}
always holds for all $k\geq 1$. 

To solve problem \cref{eq:sym_tensor_diagonalization-gener-1}, which is defined on the special linear group $\SL_{m}(\CC)$, the \emph{gradient-based Jacobi-type} (Jacobi-G) algorithms will be proposed in  \Cref{subsecc:bcd-g} (more detailedly in \Cref{subsec:Jaco_G}), which can be seen as non-orthogonal analogues of the Jacobi-G algorithms on orthogonal group \cite{IshtAV13:simax,LUC2017globally} or unitary group \cite{ULC2019}.
Then their global and weak convergence will be established in \Cref{sec:BCD-G} using the \L{}ojasiewicz gradient inequality, under the assumption that the iterates $\matr{X}_{k}$ are bounded, that is, 
there exists a universal positive constant $\mathrm{M}_{\matr{X}}>0$ such that 
\begin{equation}\label{eq:condi_bounded_X_gene-2}
\|\matr{X}_{k}\|\leq \mathrm{M}_{\matr{X}}
\end{equation}
always holds for all $k\geq 1$. 
To our knowledge, this is the first time that the theoretical convergence is established for the Jacobi-type algorithms on $\SL_{m}(\CC)$. 

\subsection{Organization}
The paper is organized as follows. In \Cref{sec:back_summary}, we present the BCD-G and Jacobi-G algorithms, define four kinds of elementary transformations and give a summary of the main results. 
In \Cref{sec:geometr_basics}, we recall the basics
of first-order geometries on the Stiefel manifold $\St(m,n,\CC)$ and special linear group $\SL_{m}(\CC)$, as well as the convergence results related to \L{}ojasiewicz inequality. 
In \Cref{subsec:first_block_algo}, we show the details of how to use the line search descent method to update the first block variable in $\St(m,n,\CC)$. 
In \Cref{sec:Jacobi_g_rota}, we define four kinds of elementary functions and present the details of three subalgorithms. 
In \Cref{sec:trian_diago_rotation} and \Cref{sec:trian_diago_rotation-2}, we present the details of four kinds of elementary transformations for JADM problem. 
In \Cref{sec:BCD-G}, 
we prove our main results about the global and weak convergence of BCD-G and Jacobi-G algorithms. 
In \Cref{sec-experiment-hooi}, some experiments are conducted to compare the proposed algorithms. 
\Cref{sec:conclusion} concludes this paper with some final remarks and possible future work. 

\section{Gradient-based algorithmic framework and a summary of results}\label{sec:back_summary}

\subsection{BCD-G and Jacobi-G algorithms}\label{subsecc:bcd-g}

Suppose that $\{\mathcal{M}_i\}_{1\leq i\leq d}$ are smooth manifolds. 
To minimize a smooth function 
\begin{equation}\label{definition-f-general}
\tilde{f}:\ \mathcal{M}_1\times\mathcal{M}_2\times\cdots\times\mathcal{M}_d \longrightarrow \RR^{+},
\end{equation}
a popular approach is the \emph{block coordinate descent} (BCD) algorithm  \cite{bertsekas1997nonlinear,luo1992convergence,luo1993error,wright2015coordinate,xu2013block,li2019polar}.
In this method, only one block variable is updated at each iteration, while other block variables are fixed; in other words, 
the problem \cref{definition-f-general} is decomposed into a sequence of lower-dimensional optimization problems. 
In the BCD algorithm, 
there are different ways to choose blocks for optimization, including the \emph{essentially cyclic, cyclic}, \emph{random} fashions \cite{wright2015coordinate,xu2013block} and the so-called 
\emph{maximum block improvement} (MBI) method \cite{chen2012maximum,li2015convergence}. 

If $d=2$, $\mathcal{M}_1=\St(m,n,\CC)$ and $\mathcal{M}_2=\SL_{m}(\CC)$, 
then problem \cref{definition-f-general} reduces to our cost function \cref{eq:sym_tensor_diagonalization-gener-0-equi}. 
For $\matr{\omega}=(\matr{Y},\matr{X})\in\St(m,n,\CC)\times\SL_{m}(\CC)$, 
we denote 
\begin{align}\label{eq:restrc_functions}
f_{1,\matr{X}}: \matr{Y}\mapsto f(\matr{Y},\matr{X}),\ \  f_{2,\matr{Y}}: \matr{X}\mapsto f(\matr{Y},\matr{X}), 
\end{align}
as the two restricted functions, which are defined on $\St(m,n,\CC)$ and  $\SL_{m}(\CC)$, respectively. 
For simplicity, we denote their  \emph{Riemannian gradients}\footnote{See \cite[Section 3.6]{absil2009optimization} and \Cref{sec:geometr_basics} for a detailed definition.} as 
$\ProjGrad{f_{1}}{\matr{\omega}}\eqdef\ProjGrad{f_{1,\matr{X}}}{\matr{Y}}$ and $\ProjGrad{f_{2}}{\matr{\omega}}\eqdef\ProjGrad{f_{2,\matr{Y}}}{\matr{X}}$,
and the Riemannian gradient of $f$ in \cref{eq:sym_tensor_diagonalization-gener-0-equi} at $\matr{\omega}$ as $\ProjGrad{f}{\matr{\omega}}$. 
To minimize the function \cref{eq:sym_tensor_diagonalization-gener-0-equi}, 
we now propose the following \emph{gradient-based block coordinate descent} (BCD-G) algorithm in \Cref{alg:BCD-G-RL}. 

\begin{algorithm}[ht!]
\caption{BCD-G algorithm}\label{alg:BCD-G-RL}
\begin{algorithmic}[1]
\STATE{{\bf Input:} A starting point $\matr{\omega}_{0}=(\matr{Y}_0,\matr{X}_0)$, a positive constant $0<\upsilon<\sqrt{2}/2$.}
\STATE{{\bf Output:} Sequence of iterates $\matr{\omega}_{k}=(\matr{Y}_k,\matr{X}_k)$.} 
\FOR{$k=1,2,\cdots,$}
\STATE Choose $t_k=1$ or $2$ such that the Riemannian gradients satisfy 
\begin{equation}
\|\ProjGrad{f_{t_k}}{\matr{\omega}_{k-1}}\|\geq\upsilon\|\ProjGrad{f}{\matr{\omega}_{k-1}}\|;
\label{eq:inequality_BCD_G}
\end{equation}
\vspace{-0.5cm}
\IF {$t_k=1$} 
\STATE Update $\matr{Y}_{k}$ using the line search descent method (cf. \Cref{subsec:line_seach});
\STATE Set $\matr{X}_{k}=\matr{X}_{k-1}$;
\ELSE
\STATE Set $\matr{Y}_{k}=\matr{Y}_{k-1}$; 
\STATE Update $\matr{X}_{k}$ using  elementary transformations (cf. \autoref{alg:jacobi-LU-G_c} to \ref{alg:jacob-GU-G_c}).
\ENDIF 
\ENDFOR
\end{algorithmic}
\end{algorithm}

In each iteration of \Cref{alg:BCD-G-RL}, instead of the frequently used cyclic or random fashion to choose the block for optimization, we choose the block $t_k=1$ or $2$ satisfying the inequality\footnote{The inequality \cref{eq:inequality_BCD_G} can be seen as a block coordinate analogue of \cite[Eq. (3.3)]{IshtAV13:simax} and \cite[Eq. (10)]{LUC2017globally}.} \cref{eq:inequality_BCD_G}. 
Since the Riemannian gradients are related as 
\begin{align}\label{eq:grad_relation}
\ProjGrad{f}{\matr{\omega}} = (\ProjGrad{f_1}{\matr{\omega}},\ \ProjGrad{f_2}{\matr{\omega}}),
\end{align}
we have that 
$\|\ProjGrad{f}{\matr{\omega}}\|^2=\|\ProjGrad{f_1}{\matr{\omega}}\|^2+\|\ProjGrad{f_2}{\matr{\omega}}\|^2.$ 
Therefore, in each iteration, if $0<\upsilon<\sqrt{2}/2$, we can always choose $t_k=1$ or $2$ such that the inequality \cref{eq:inequality_BCD_G} is satisfied, and thus \Cref{alg:BCD-G-RL} is well defined.

In \Cref{alg:BCD-G-RL}, to update $\matr{Y}_k$, we choose the \emph{line search descent} method \cite{absil2005convergence,absil2009optimization,nocedal2006numerical,ring2012optimization,sato2015new}, which will be detailedly presented in \Cref{subsec:first_block_algo}. 
To update $\matr{X}_k$, as in Jacobi-type methods, we use four kinds of elementary transformations (will be detailed introduced in \Cref{subsec:elementary_trans}), including the \emph{Givens plane, plane upper triangular, plane lower triangular} and \emph{plane diagonal transformations}\footnote{The reason why we use plane diagonal transformations will be shown in \Cref{sec:Jacobi_g_rota}.}. 
We group these elementary transformations into three classes (GLU, GQU and GU)  motivated by well-known matrix decompositions,
which give rise to three different variants of  \Cref{alg:BCD-G-RL} (BCD-GLU, BCD-GQU and BCD-GU).
We recall the matrix decompositions and related Lie groups in \Cref{subsec:matr_dec},
before introducing the elementary transformations and their classes in \Cref{subsec:elementary_trans}.

Similarly to \Cref{alg:BCD-G-RL}, we propose optimization algorithms for minimization of the cost function \cref{eq:sym_tensor_diagonalization-gener-1} for the square case (second reformulation on $\SL_{m}(\CC)$).
In these algorithms, $\matr{X}_k$ is updated with  four elementary transformations, and therefore they are Jacobi-type algorithms.
We summarize these \emph{gradient-based Jacobi-type} (Jacobi-G) algorithms in \Cref{alg:jacobi-LU-G_c-1}.

\begin{algorithm}
\caption{Jacobi-G algorithm}\label{alg:jacobi-LU-G_c-1}
\begin{algorithmic}[1]
\STATE{{\bf Input:} A  starting point $\matr{X}_{0}$.}
\STATE{{\bf Output:} Sequence of iterates $\{\matr{X}_{k}\}_{k\ge1}$.}
\FOR{$k=1,2,\cdots,$}
\STATE  Update $\matr{X}_k$ using elementary transformations (cf. \autoref{alg:jacobi-LU-G_c} to \ref{alg:jacobi-LQ-G_c}).
\ENDFOR
\end{algorithmic}
\end{algorithm}

\Cref{alg:jacobi-LU-G_c-1} can be seen as a non-orthogonal analogue of the Jacobi-G algorithm in \cite{IshtAV13:simax,LUC2017globally,ULC2019}.
As with BCD-G, two types of Jacobi-G exist: Jacobi-GLU and Jacobi-GQU, based on GLU and GQU classes of elementary transformations, respectively.
Roughly speaking, these  algorithms are variants of \Cref{alg:BCD-G-RL}, where only  $\matr{X}_k$ is updated.

\subsection{Matrix decompositions and matrix groups}\label{subsec:matr_dec}
A matrix $\matr{X}\in\CC^{m\times m}$ is said to be \emph{upper triangular} if $\matrelem{X}_{ij}=0$ for $i>j$.
Let $\UT_{m}(\CC)\subseteq\GL_{m}(\CC)$ be the \emph{upper triangular subgroup}. 
Let $\EUT_{m}(\CC)=\UT_{m}(\CC)\cap\SL_{m}(\CC)$, \emph{i.e.}, the set of upper triangular matrices with determinant equal to 1. 
Similarly, we 
let $\LT_{m}(\CC)\subseteq\CC^{m\times m}$ be the \emph{lower triangular subgroup} and  $\ELT_{m}(\CC)=\LT_{m}(\CC)\cap\SL_{m}(\CC)$.
Let $\UNN{m}(\CC)\subseteq \CC^{m\times m}$ be the \emph{unitary group}, and $\SUNN{m}(\CC)\subseteq\UNN{m}(\CC)$ be the \emph{special unitary group}. 

We first discuss the matrix decompositions of $\SL_{m}(\CC)$. 

\begin{itemize}
\item Any matrix $\matr{X}\in\SL_{m}(\CC)$ has the \emph{LU decomposition} \cite{GoluV96:jhu} $\matr{X}=\matr{L}\matr{U}$ with $\matr{L}\in\LT_{m}(\CC)$ and $\matr{U}\in\UT_{m}(\CC)$. 
We use the shorthand notation 
\begin{align}
\SL_m(\CC) =
\ELT_{m}(\CC)\bullet\EUT_{m}(\CC),
\label{eq:rela-equiavalent-5}
\end{align}
where $\mathcal{A} \bullet \mathcal{B}$ denotes the set of all matrix product for matrices coming from two matrix sets $\mathcal{A}$ and $\mathcal{B}$.
The decomposition \cref{eq:rela-equiavalent-5} motivates the GLU class, which includes the plane lower triangular, plane upper triangular and plane diagonal transformations (\autoref{alg:jacobi-LU-G_c}), and is used in BCD-GLU and Jacobi-GLU algorithms.  

\item Any matrix $\matr{X}\in\SL_{m}(\CC)$ has the \emph{QU decomposition}\footnote{This is also called QR decomposition in the literature.} $\matr{X}=\matr{Q}\matr{U}$ with $\matr{Q}\in\SUNN{m}(\CC)$ and $\matr{U}\in\UT_{m}(\CC)$, which can be compactly written as 
\begin{align}
\SL_m(\CC) =
\SUNN{m}(\CC)\bullet\EUT_{m}(\CC).
\label{eq:rela-equiavalent-6}
\end{align}
The decomposition \cref{eq:rela-equiavalent-6} motivates the GQU class, which includes the Givens plane, plane upper triangular and plane diagonal transformations (\autoref{alg:jacobi-LQ-G_c}), and is used in BCD-GQU and Jacobi-GQU algorithms. 
\end{itemize}

The decompositions mentioned above  can be used to parameterize $\textbf{RSL}(m,n,\CC)$.
Indeed,  \Cref{lem:rela-equiavalent} in the compact notation can be written as
\[
\textbf{RSL}(m,n,\CC) =\St(m,n,\CC)\bullet\SL_m(\CC),
\]
which gives rise to  LU- and QU-based decompositions of $\textbf{RSL}(m,n,\CC)$:
\begin{align}
\textbf{RSL}(m,n,\CC) &=
\St(m,n,\CC)\bullet\ELT_{m}(\CC)\bullet\EUT_{m}(\CC),
\label{eq:rela-equiavalent-2} \\
\textbf{RSL}(m,n,\CC) &=
\St(m,n,\CC)\bullet\SUNN{m}(\CC)\bullet\EUT_{m}(\CC).
\label{eq:rela-equiavalent-3}
\end{align}
Moreover, for  $\textbf{RSL}(m,n,\CC)$, a third decomposition is possible, using the fact that  $\St(m,n,\CC)\bullet\SUNN{m}(\CC) = \St(m,n,\CC)$.
Then the equation \cref{eq:rela-equiavalent-3} can be simplified as 
\begin{align}
\textbf{RSL}(m,n,\CC) =\St(m,n,\CC)\bullet\EUT_{m}(\CC),
\label{eq:rela-equiavalent-13}
\end{align}
which can also be interpreted as applying the  QU decomposition  to a rectangular matrix from  $\textbf{RSL}(m,n,\CC)$. 
This gives rise to the third class GU, which only includes the plane upper triangular and plane diagonal transformations (\autoref{alg:jacob-GU-G_c}), and is used in BCD-GU. 
%In this case, we call \Cref{alg:BCD-G-RL} the BCD-GU algorithm. 

\subsection{Elementary transformations}\label{subsec:elementary_trans}
Let us introduce a few more matrix groups. 
An upper triangular matrix $\matr{X}$ is said to be \emph{unipotent} if it satisfies $\matrelem{X}_{ii}=1$ for $1\leq i\leq m$. 
Let $\SUT_{m}(\CC)\subseteq \EUT_{m}(\CC) $ be the \emph{upper unipotent subgroup} of unipotent upper triangular matrices. 
Similarly,
we let $\SLT_{m}(\CC)\subseteq \ELT_{m}(\CC)$ be the \emph{lower unipotent subgroup}.
Finally, a diagonal matrix $\matr{X}\in\CC^{m\times m}$ is said to be a \emph{diagonal transformation} if the product of all the diagonal elements is equal to 1.
Let $\DT_{m}(\CC)\subseteq\GL_{m}(\CC)$ be the set of diagonal transformation matrices.

The elementary transformations are based on the following $2\times 2$ matrices:   
{\small
\begin{align*}
\SUT_{2}(\CC) = 
\left\{\begin{bmatrix}
1 &  z\\
0& 1
\end{bmatrix},
z\in\CC\right\},\ 
\SLT_{2}(\CC) = 
\left\{\begin{bmatrix}
1 &  0\\
z& 1
\end{bmatrix},
z\in\CC\right\},\  
 \DT_{2}(\CC) = 
\left\{\begin{bmatrix}
z &  0\\
0& \frac{1}{z}
\end{bmatrix},
z\in\CC_{*}\right\},
\end{align*}}
as well as the $2\times 2$ matrices from $\SUNN{2}(\CC)$.

Let $(i,j)$ be a pair of indices satisfying $1 \le i < j \le m$. 
We introduce an operator $\mathcal{E}_{i,j}: \CC^{2\times 2} \to\CC^{m\times m}$ sending $\Psi$ to $\matr{X} \in \CC^{m\times m}$ satisfying 
%  \begin{align*}
% &\matrelem{X}_{\ell\ell}  = 1, \quad \ell \not\in \{i,j\}, \\
% &\begin{bmatrix}
% \matrelem{X}_{ii} &  \matrelem{X}_{ij}\\
% \matrelem{X}_{ji}& \matrelem{X}_{jj}
% \end{bmatrix}  = \begin{bmatrix}
% \matrelem{\Psi}_{11} &  \matrelem{\Psi}_{12}\\
% \matrelem{\Psi}_{21}& \matrelem{\Psi}_{22}
% \end{bmatrix}, \\
% &\matrelem{X}_{k\ell} = 0, \quad \text{otherwise}. 
% \end{align*}
 \begin{align*}
\begin{bmatrix}
\matrelem{X}_{ii} &  \matrelem{X}_{ij}\\
\matrelem{X}_{ji}& \matrelem{X}_{jj}
\end{bmatrix}  = \begin{bmatrix}
\matrelem{\Psi}_{11} &  \matrelem{\Psi}_{12}\\
\matrelem{\Psi}_{21}& \matrelem{\Psi}_{22}
\end{bmatrix}, \ \
\bigg\{
\begin{aligned}
\matrelem{X}_{\ell\ell}  = 1,\ \ & \text{if}\ \ell \not\in \{i,j\}, \\
\matrelem{X}_{k\ell} = 0,\ \ & \text{otherwise}. 
\end{aligned}
\end{align*} 
%$\mathcal{E}_{i,j}(\matr{\Utwo})$ 
% \begin{equation*}{\small
% \mathcal{E}_{i,j}(\matr{\Utwo}) \eqdef
% \begin{matrix}
% & & & i & & j & & & \\
% \ldelim[{8}{0.1cm} & 1 & & & & & & & \rdelim]{8}{0.1cm} & \\
% & & \ddots & & & & \mathbf{0} & & & \\
% & & & \Psi_{11} & & \Psi_{12} & & & & i \\
% & & &  & \ddots & & & & & \\
% & & & \Psi_{21} & & \Psi_{22} & & & & j\\
% & & \mathbf{0} & & & & \ddots & & & \\
% & & & & & & & 1 & & \\
% \end{matrix}}
% \end{equation*}
% is the identity matrix $\matr{I}_m$ except that $\mathcal{P}_{i,j}(\mathcal{E}_{i,j}(\matr{\Utwo}))=\matr{\Utwo}$. 
Now we define the following four elementary transformations on $\SL_{m}(\CC)$: 
\begin{itemize}
\item $\Gmat{i}{j}{\matr{\Utwo}}\eqdef\mathcal{E}_{i,j}(\matr{\Utwo})$: \emph{Givens plane transformation} for $\matr{\Utwo} \in \SUNN{2}(\CC)$;
\item $\Umat{i}{j}{\matr{\Utwo}}\eqdef\mathcal{E}_{i,j}(\matr{\Utwo})$: \emph{plane upper triangular transformation} for $\matr{\Utwo} \in \SUT_{2}(\CC)$;
\item $\Lmat{i}{j}{\matr{\Utwo}}\eqdef\mathcal{E}_{i,j}(\matr{\Utwo})$: \emph{plane lower triangular transformation} for $\matr{\Utwo} \in \SLT_{2}(\CC)$;
\item $\Dmat{i}{j}{\matr{\Utwo}}\eqdef\mathcal{E}_{i,j}(\matr{\Utwo})$: \emph{plane diagonal transformation} for $\matr{\Utwo} \in \DT_{2}(\CC)$.
\end{itemize}

\begin{remark}
These elementary transformations have all been used in the literature. The Givens  transformations $\Gmat{i}{j}{\matr{\Utwo}}$ were used very often in the Jacobi-type algorithms for joint approximate diagonalization of matrices or tensors by orthogonal or non-orthogonal transformations  \cite{Como10:book,LUC2017globally,ULC2019,Afsa07:honolulu,afsari2006simple,sorensen2009approximate}. 
Triangular transformations $\Umat{i}{j}{\matr{\Utwo}}$ and  $\Lmat{i}{j}{\matr{\Utwo}}$ also appeared  many times in the Jacobi-type algorithms on $\SL_{m}(\CC)$ or $\SL_{m}(\RR)$ \cite{afsari2004gradient,afsari2006simple,maurandi2013fast,maurandi2014decoupled,maurandi2014jacobi}. 
In the real case, the diagonal transformation $\Dmat{i}{j}{\matr{\Utwo}}$ was once used in \cite{souloumiac2009nonorthogonal}.
\end{remark}

The iterates $\matr{X}_k$ in \Cref{alg:BCD-G-RL} and \Cref{alg:jacobi-LU-G_c-1} are updated multiplicatively as $\matr{X}_k = \matr{X}_{k-1}\matr{P}_{k}$, where $\matr{P}_{k}$ is an elementary transformation for a pair of indices $(i_k,j_k)$ belonging to one of the following three classes. 
These three classes are inspired by equations \cref{eq:rela-equiavalent-2}, \cref{eq:rela-equiavalent-3}, \cref{eq:rela-equiavalent-13} and by a similar idea as in \cite{IshtAV13:simax,LUC2017globally,ULC2019}.
We call them the \emph{GLU} (based on LU decomposition), \emph{GQU} (based on QU decomposition) and GU transformations, respectively. 
\begin{itemize}
\item {GLU class}:  $\matr{P}_{k}=\Lmat{i_k}{j_k}{\matr{\Utwo}^{*}_{k}},\ \Umat{i_k}{j_k}{\matr{\Utwo}^{*}_{k}}$ or $\Dmat{i_k}{j_k}{\matr{\Utwo}^{*}_{k}}$;
\item GQU class:  $\matr{P}_{k}=\Gmat{i_k}{j_k}{\matr{\Utwo}^{*}_{k}},\ \Umat{i_k}{j_k}{\matr{\Utwo}^{*}_{k}}$ or $\Dmat{i_k}{j_k}{\matr{\Utwo}^{*}_{k}}$;
\item GU class:  $\matr{P}_{k}= \Umat{i_k}{j_k}{\matr{\Utwo}^{*}_{k}}$ or $\Dmat{i_k}{j_k}{\matr{\Utwo}^{*}_{k}}$.
\end{itemize}
The choice of the pair $(i_k,j_k)$, the matrix $\matr{\Utwo}^{*}_{k}$ and the particular type of transformations in each class will be given in  \autoref{alg:jacobi-LU-G_c}, \autoref{alg:jacobi-LQ-G_c} and \autoref{alg:jacob-GU-G_c}.
The algorithms and their convergence results are summarized in  \Cref{table-example-3-0}. 

{\renewcommand{\arraystretch}{1.2}

\begin{table}[h!]
\centering
\caption{A summary of the proposed algorithms}
\label{table-example-3-0}
\scalebox{0.8}{
\begin{tabular}{| p{4.0cm} | p{2.0cm} | p{2.7cm} | p{2.2cm} | p{2.0cm} | p{2.0cm} |}
\toprule
$\makecell[c]{\textrm{Model}}$ & $\textrm{\makecell[l]{Proposed\\ algorithms}}$  &
$\makecell[c]{\textrm{Location}}$  &$\textrm{\makecell[l]{Elementary\\ transformations}}$  & $\textrm{\makecell[l]{Global\\ convergence}}$ & $\textrm{\makecell[l]{Weak\\ convergence}}$\\
\hline
\multirow{3}{*}{\textrm{\makecell[c]{First reformulation  \cref{eq:sym_tensor_diagonalization-gener-0-equi}\\ on $\St(m,n,\CC)\times\SL_{m}(\CC)$}}} & \textrm{BCD-GLU} & \Cref{alg:BCD-G-RL} \& \autoref{alg:jacobi-LU-G_c} & $\matr{L}, \matr{U}, \matr{D}$ & \multirow{3}{*}{\Cref{theorem_main_global_conv}}&\multirow{3}{*}{\Cref{theorem_main_weak_conv}}\\
\cline{2-4}
& \textrm{BCD-GQU} & \Cref{alg:BCD-G-RL} \& \autoref{alg:jacobi-LQ-G_c} & $\matr{Q}, \matr{U}, \matr{D}$& &\\
\cline{2-4}
& \textrm{BCD-GU} & \Cref{alg:BCD-G-RL} \& \autoref{alg:jacob-GU-G_c} & $\matr{U}, \matr{D}$&&\\
\hline
\multirow{2}{*}{\textrm{\makecell[c]{Second reformulation  \cref{eq:sym_tensor_diagonalization-gener-1}\\ on $\SL_{m}(\CC)$}}}
& \textrm{Jacobi-GLU} & \Cref{alg:jacobi-LU-G_c-1} \& \autoref{alg:jacobi-LU-G_c} & $\matr{L}, \matr{U}, \matr{D}$ & \multirow{2}{*}{\Cref{theorem:weak_Jaco_Lu-0}} & \multirow{2}{*}{\Cref{theorem_main_weak_conv-Jacobi-GLU-09}}\\
\cline{2-4}
& \textrm{Jacobi-GQU} &  \Cref{alg:jacobi-LU-G_c-1} \& \autoref{alg:jacobi-LQ-G_c} & $\matr{Q}, \matr{U}, \matr{D}$&&\\
\bottomrule
\end{tabular}}
\end{table}}

\begin{remark}
While the algorithms and convergence results described in this paper are provided for complex matrices, complex Stiefel manifold $\St(m,n,\CC)$ and complex special linear group $\SL_{m}(\CC)$, they also remain valid in the real case. 
\end{remark}

\section{Geometries on $\St(m,n,\CC)$ and $\SL_{m}(\CC)$}\label{sec:geometr_basics}

\subsection{Notations} 
Let $1\leq m\leq n$. 
For a complex matrix $\matr{Z}\in \CC^{n\times m}$ and a complex number $z\in \CC$, we write the real and imaginary parts as $\matr{Z} =\matr{Z}^{\R} + \ui \matr{Z}^{\I}$ and $z=\R(z)+\ui\I(z)$, respectively. 
For complex matrices $\matr{Z}_1,\matr{Z}_2\in \CC^{n\times m}$, we introduce the following real-valued inner product
\begin{equation}\label{eq:RInner}
\RInner{\matr{Z}_1}{\matr{Z}_2} \eqdef 
\langle\matr{Z}_1^{\R},\matr{Z}_2^{\R} \rangle + \langle\matr{Z}_1^{\I},\matr{Z}_2^{\I} \rangle  =
\Re \left({\trace{\matr{Z}_1^{\H}\matr{Z}_2}}\right),
\end{equation}
which makes $\CC^{n\times m}$ a real Euclidean space  of dimension $2nm$.
Let $h:\CC^{n\times m}\rightarrow\RR$ be a differentiable function and $\matr{Z}\in\CC^{n\times m}$. 
We denote by  $\frac{\partial h}{\partial\matr{Z}^\R},\frac{\partial h}{\partial\matr{Z}^\I} \in \RR^{n\times m}$ the matrix Euclidean derivatives of $h$ with respect to real and imaginary parts of $\matr{Z}$. 
The \emph{Wirtinger derivatives} \cite{abrudan2008steepest,brandwood1983complex,krantz2001function} are defined as
\begin{equation*}
\frac{\partial h}{\partial\matr{Z}^{*}} \eqdef \frac{1}{2}\left(\frac{\partial h}{\partial\matr{Z}^\R}+ \ui\frac{\partial h}{\partial\matr{Z}^\I}\right), \quad 
\frac{\partial h}{\partial\matr{Z}} \eqdef \frac{1}{2}\left(\frac{\partial h}{\partial\matr{Z}^\R}- \ui\frac{\partial h}{\partial\matr{Z}^\I}\right). 
\end{equation*}
Then the Euclidean gradient of $h$ with respect to the inner product \cref{eq:RInner} becomes
\begin{equation}\label{eq:matr_Euci_grad}
\nabla h(\matr{Z}) = \frac{\partial h}{\partial\matr{Z}^\R}+ \ui\frac{\partial h}{\partial\matr{Z}^\I} = 2 \frac{\partial h}{\partial\matr{Z}^{*}}.
\end{equation}
For real matrices $\matr{Z}_1,\matr{Z}_2\in\RR^{n\times m}$, we see that \cref{eq:RInner} becomes the standard inner product, and \cref{eq:matr_Euci_grad} becomes the standard Euclidean gradient. 
We denote by $\mathbb{S}_{2}\subseteq\RR^{3}$ the unit sphere, and $\CC_{*}=\CC\backslash \{0\}$.

\subsection{Riemannian gradient on $\St(m,n,\CC)$}

For a matrix  $\matr{C}\in\CC^{m\times m}$, we denote 
$\symmm{\matr{C}}\eqdef\frac{1}{2}(\matr{C}+\matr{C}^{\H})$ 
and 
$\skeww{\matr{C}}\eqdef\frac{1}{2}(\matr{C}-\matr{C}^{\H})$. 
Let ${\rm\bf T}_{\matr{Y}} \St(m,n,\CC)$ be the \emph{tangent space} to $\St(m,n,\CC)$ at a point $\matr{Y}\in\St(m,n,\CC)$.  
Let $\matr{Y}_{\perp}\in\CC^{n\times (n-m)}$ be an orthogonal complement of $\matr{Y}$, that is, $[\matr{Y},  \matr{Y}_{\perp}]\in\CC^{n\times n}$ is a unitary matrix.  
By \cite[Definition 6]{manton2001modified}, we know that 
\begin{align*}{\small
{\rm\bf T}_{\matr{Y}} \St(m,n,\CC) 
= \{\matr{V}\in\CC^{n\times m}, \matr{V}=\matr{Y}\matr{C}+\matr{Y}_{\perp}\matr{B},
\matr{C}\in\CC^{m\times m}, \matr{C}^{\H}+\matr{C}=0, \matr{B}\in\CC^{(n-m)\times m}\},\label{eq:tanget_spac}}
\end{align*} 
which is a $(2nm-m^2)$-dimensional vector space. 
The orthogonal projection of a matrix  $\xi\in\CC^{n\times m}$ onto ${\rm\bf T}_{\matr{Y}} \St(m,n,\CC)$ is 
\begin{equation}\label{eq:proj_Stiefel}
{\rm Proj}_{\matr{Y}} \xi = (\matr{I}_{n}-\matr{Y}\matr{Y}^{\H})\xi + \matr{Y}\skeww{\matr{Y}^{\H}\xi} 
=\xi - \matr{Y}\symmm{\matr{Y}^{\H}\xi}.
\end{equation} 
We denote ${\rm Proj}^{\bot}_{\matr{Y}} \xi \eqdef\xi - {\rm Proj}_{\matr{Y}} \xi $. 
Let $p:\St(m,n,\CC)\rightarrow\RR$ be a differentiable function, and $\matr{Y}\in\St(m,n,\CC)$. 
Note that $\St(m,n,\CC)$ is an embedded submanifold of the Euclidean space $\CC^{n\times m}$. 
By equation \cref{eq:proj_Stiefel}, we have the Riemannian gradient of $p$ at $\matr{Y}$ as:
\begin{align}\label{eq:Rie_grad}
\ProjGrad{p}{\matr{Y}} = {\rm Proj}_{\matr{Y}} \nabla p(\matr{Y})
= \nabla p(\matr{Y}) -\matr{Y}\symmm{\matr{Y}^{\H}\nabla p(\matr{Y})}.
\end{align}
By \cite[Example 5.4.2]{absil2009optimization},
the \emph{exponential map} at $\matr{Y}$ is defined as
\begin{align}\label{eq:ex_sl_pr}{\small
{\rm Exp}_{\matr{Y}}: {\rm\bf T}_{\matr{Y}} \St(m,n,\CC)}&{\small\longrightarrow \St(m,n,\CC)}\\
{\small\matr{V}}&{\small\longmapsto [\matr{Y}, \matr{V}]
\expp{
\begin{bmatrix}
\matr{Y}^{H}\matr{V} &-\matr{V}^{\H}\matr{V}\\
\matr{I}_{m} &\matr{Y}^{H}\matr{V}
\end{bmatrix}}
\begin{bmatrix}
\expp{-\matr{Y}^{H}\matr{V}} \\
\matr{0}_{m\times m}
\end{bmatrix},\notag}
\end{align}
where 
$\exp(\cdot)$ is the matrix exponential function \cite{absil2009optimization,baker2012matrix,Hall15:lie}.

\subsection{Riemannian gradient on $\SL_{m}(\CC)$}
\label{subsec:grad_SL_n}
Let $\sll_{m}(\CC)\eqdef\{\matr{X}\in\CC^{m\times m}, \tr(\matr{X})=0\}$ be the \emph{Lie algebra} \cite{baker2012matrix} of the complex special linear group $\SL_{m}(\CC)$. 
Then the tangent space to $\SL_{m}(\CC)$ at a point $\matr{X}\in\SL_{m}(\CC)$ can be constructed \cite[Eq. (3.7),(3.8)]{baker2012matrix} by 
\begin{equation}\label{eq:tang_slm}
\mathbf{T}_{\matr{X}}\SL_{m}(\CC) = \{\matr{X}\Omega, \Omega\in\sll_{m}(\CC)\}.
\end{equation} 
Let
$\su_{m}(\CC)\eqdef\{\matr{X}\in\CC^{m\times m}, \matr{X}^{\H} = -\matr{X}, \tr(\matr{X})=0\}$ be the Lie algebra of the special unitary group $\SUNN{m}(\CC)$.  
Then the tangent space to $\SUNN{m}(\CC)$ at a point $\matr{X}\in\SUNN{m}(\CC)$ can be constructed \cite[Eq. (3.15)]{baker2012matrix} by 
$\mathbf{T}_{\matr{X}}\SUNN{m}(\CC) = \{\matr{X}\Omega, \Omega\in\sll_{m}(\CC)\}.$

Let $\mathbf{T}_{\matr{X}}\SL_{m}(\CC)$ be the tangent space to $\SL_{m}(\CC)$ at a point $\matr{X}\in\SL_{m}(\CC)$ as in \cref{eq:tang_slm}. 
For tangent matrices $\matr{V}_1, \matr{V}_2\in\mathbf{T}_{\matr{X}}\SL_{m}(\CC)$,
we use the \emph{left invariant}  \cite{absil2009optimization}, \cite[Eq. (6.2)]{afsari2004gradient}  Riemannian metric  
\begin{equation*}
\langle\matr{V}_1, \matr{V}_2\rangle_{\matr{X}} 
\eqdef \RInner{\matr{X}^{-1}\matr{V}_1}{\matr{X}^{-1}\matr{V}_2}
= \Re \left({\trace{\matr{V}_1^{\H}(\matr{X}\matr{X}^{\H})^{-1}\matr{V}_2}}\right). 
\end{equation*}
Let
$\textit{g}: \SL_{m}(\CC) \longrightarrow \RR^{+}$  
be a differentiable function, and $\matr{X}\in\SL_{m}(\CC)$. 
Then the Riemannian gradient of $g$ at $\matr{X}$ is the orthogonal projection \cite[Lemma 6.2]{afsari2004gradient} of its Euclidean gradient $\nabla g(\matr{X})$ to $\mathbf{T}_{\matr{X}}\SL_{m}(\CC)$, that is,
\begin{align}\label{eq:riema_gradient_c}
\ProjGrad{g}{\matr{X}} = \matr{X}\left(\matr{X}^{\H}\nabla g(\matr{X}) - \frac{\tr(\matr{X}^{\H}\nabla g(\matr{X}))}{n}\matr{I}_n\right).
\end{align}
We denote $\matr{\Lambda}(\matr{X})\eqdef\matr{X}^{-1}\ProjGrad{g}{\matr{X}}\in\sll_{m}(\CC)$ for $\matr{X}\in\SL_{m}(\CC)$, which will be frequently used in this paper. 

In what follows, we will use the following exponential map  
\begin{equation}\label{eq:Exp_x_sl}
\text{Exp}_{\matr{X}}: \mathbf{T}_{\matr{X}}\SL_{m}(\CC)  \to \SL_{m}(\CC),\ \ \matr{X}\Omega\mapsto\matr{X} \exp(\Omega), 
\end{equation}
where 
$\exp(\cdot)$ is the matrix exponential function \cite{absil2009optimization,baker2012matrix,Hall15:lie}.
For any tangent matrix $\matr{V}\in \mathbf{T}_{\matr{X}}\SL_{m}(\CC)$, we have the following relationship between $\text{Exp}_{\matr{X}}$ in \cref{eq:Exp_x_sl} and the Riemannian gradient \cite[Eq. (3.31)]{absil2009optimization}: 
\begin{equation}\label{eq:gradient_exp_map}
\langle\matr{V},\ProjGrad{g}{\matr{X}}\rangle_{\matr{X}}  =
\left.\left( \frac{d}{dt} g (\text{Exp}_{\matr{X}}(t\matr{V}) )\right)\right|_{t=0},
\end{equation}
which will be used in the proof of \cref{lem:ProjGradSubmatrix}. 

\subsection{Tangent spaces to other matrix groups}\label{subsec:tang_22}

A matrix $\matr{X}\in\CC^{m\times m}$ is said to be \emph{strictly upper triangular} if $\matrelem{X}_{ij}=0$ for $i\geq j$.
Let $\sut_{m}(\CC)\subseteq\CC^{m\times m}$ be the set of strictly upper triangular matrices.
Then the tangent space to $\SUT_{m}(\CC)$ at a point $\matr{X}\in\SUT_{m}(\CC)$ can be constructed \cite[Eq. (3.11)]{baker2012matrix}, \cite[Section 6.4]{afsari2004gradient} by $\mathbf{T}_{\matr{X}}\SUT_{m}(\CC) = \{\matr{X}\Omega, \Omega\in\sut_{m}(\CC)\}$. 
Similar as above, we let  $\slt_{n}(\CC)\subseteq\CC^{m\times m}$ be the set of \emph{strictly lower triangular} matrices.
Then the tangent space to $\SLT_{m}(\CC)$ at a point $\matr{X}\in\SLT_{m}(\CC)$ can be constructed by $\mathbf{T}_{\matr{X}}\SLT_{m}(\CC) = \{\matr{X}\Omega, \Omega\in\slt_{m}(\CC)\}$.
Let $\sdt_{m}(\CC)\subseteq\CC^{m\times m}$ be the set of diagonal traceless matrices. 
Then the tangent space to $\DT_{m}(\CC)$ at a point $\matr{X}\in\DT_{m}(\CC)$ can be constructed by $\mathbf{T}_{\matr{X}}\DT_{m}(\CC) = \{\matr{X}\Omega, \Omega\in\sdt_{m}(\CC)\}$. 
In particular, for the case $m=2$, we have 
{\small
\begin{align*}
\sut_{2}(\CC) = 
\left\{\begin{bmatrix}
0 &  z\\
0& 0
\end{bmatrix},
z\in\CC\right\},\ \ 
\slt_{2}(\CC) = 
\left\{\begin{bmatrix}
0 &  0\\
z& 0
\end{bmatrix},
z\in\CC\right\},\ \ 
\sdt_{2}(\CC) = 
\left\{\begin{bmatrix}
z &  0\\
0& -z
\end{bmatrix},
z\in\CC \right\}.
\end{align*}}

\subsection{Inequalities for convergence analysis}\label{sec:inequalities}
We recall some definitions and results about the \L{}ojasiewicz gradient inequality \cite{loja1965ensembles,lojasiewicz1993geometrie,absil2005convergence,Usch15:pjo}. 
These results were used in \cite{LUC2018,ULC2019} to prove the global convergence of Jacobi-G algorithms on the orthogonal and unitary groups, and will be used in this paper as well. 

\begin{definition} [{\cite[Definition 2.1]{SU15:pro}}]\label{def:Lojasiewicz}
Let $\mathcal{M} \subseteq \RR^d$ be a Riemannian submanifold,
and $\varphi: \mathcal{M} \to \RR$ be a differentiable function.
The function $\varphi: \mathcal{M} \to \RR$ is said to satisfy a \emph{\L{}ojasiewicz gradient inequality} at $\vect{x} \in \mathcal{M}$, if there exist
$\delta>0$, $\zeta\in (0,\frac{1}{2}]$ and a neighborhood $\mathcal{U}$ in $\mathcal{M}$ of $\vect{x}$ such that for all $\vect{y}\in\mathcal{U}$, it follows that 
\begin{equation}\label{eq:Lojasiewicz}
|{\varphi}(\vect{y})-{\varphi}(\vect{x})|^{1-\zeta}\leq \delta\|\ProjGrad{\varphi}{\vect{y}}\|.
\end{equation}
\end{definition}

\begin{lemma}[{\cite[Proposition 2.2]{SU15:pro}}]\label{lemma-SU15}
Let $\mathcal{M}\subseteq\RR^d$ be an analytic submanifold\footnote{See {\cite[Definition 2.7.1]{krantz2002primer}} or \cite[Definition 5.1]{LUC2018} for a definition of an analytic submanifold.} and $\varphi: \mathcal{M} \to \RR$ be a real analytic function.
Then $\varphi$ satisfies a \L{}ojasiewicz gradient inequality \cref{eq:Lojasiewicz} at any $\vect{x}\in \mathcal{M}$. 	
\end{lemma}

\begin{theorem}[{\cite[Theorem  2.3]{SU15:pro}}]\label{theorem-SU15}
Let $\mathcal{M}\subseteq\RR^d$ be an analytic submanifold
and 
$\{\vect{x}_k\}_{k\geq 1}\subseteq\mathcal{M}$.
Suppose that $\varphi$ is real analytic and, for large enough $k$,\\
(i) there exists $\sigma>0$ such that
\begin{equation*}
\varphi(\vect{x}_{k})-{\varphi}(\vect{x}_{k+1})\geq \sigma\|\ProjGrad{\varphi}{\vect{x}_k}\|\|\vect{x}_{k+1}-\vect{x}_{k}\|;
\end{equation*}
(ii) $\ProjGrad{\varphi}{\vect{x}_k}=0$ implies that $\vect{x}_{k+1}=\vect{x}_{k}$.\\
Then, if $\vect{x}_*$ is an accumulation point of $\{\vect{x}_k\}_{k\geq 1}$, it is the limit point. 
\end{theorem}

Since the special linear group $\SL_{m}(\CC)$ is not compact, the iterates $\{\matr{\omega}_{k}\}_{k\ge1}$ in \Cref{alg:BCD-G-RL} for cost function \cref{eq:sym_tensor_diagonalization-gener-0-equi} may have no accumulation point.
However, if there exists an accumulation point, we have the following result about its global convergence,
which is a direct consequence of \Cref{theorem-SU15} and inequality \cref{eq:inequality_BCD_G}. 

\begin{lemma}\label{theorem-SU15-2}
Suppose that, in \Cref{alg:BCD-G-RL} for cost function \cref{eq:sym_tensor_diagonalization-gener-0-equi}, the iterates $\{\matr{\omega}_k\}_{k\geq 1}$ satisfy that, for large enough $k$,\\
(i) there exists $\sigma>0$ such that
\begin{equation}\label{eq:sufficient_descent-2}
f(\matr{\omega}_{k-1})-f(\matr{\omega}_{k})\geq\sigma\|\ProjGrad{f_{t_k}}{\matr{\omega}_{k-1}}\|\|\matr{\omega}_{k}-\matr{\omega}_{k-1}\|;
\end{equation}
(ii) $\ProjGrad{f_{t_k}}{\matr{\omega}_{k-1}}=0$ implies that $\matr{\omega}_{k}=\matr{\omega}_{k-1}$.\\
Then, if $\matr{\omega}_*$ is an accumulation point of the iterates $\{\matr{\omega}_k\}_{k\geq 1}$, it is the limit point.
\end{lemma}

We also have the following result about its weak convergence,
which can be proved easily by inequality \cref{eq:inequality_BCD_G} and the fact that $f(\omega)\geq 0$. 

\begin{lemma}\label{theorem-weak-general_c}
In \Cref{alg:BCD-G-RL} for cost function \cref{eq:sym_tensor_diagonalization-gener-0-equi}, 
if there exists $\eta>0$ such that
\begin{equation}\label{equation-condition-weak_c}
f(\matr{\omega}_{k-1})-f(\matr{\omega}_{k})\geq\eta\|\ProjGrad{f_{t_k}}{\matr{\omega}_{k-1}}\|^2
\end{equation}
always holds, then $\lim_{k\rightarrow\infty}\ProjGrad{f}{\matr{\omega}_{k-1}}=0$.
In particular, if $\matr{\omega}_*$ is an accumulation point of the iterates $\{\matr{\omega}_{k}\}_{k\ge1}$,
then $\matr{\omega}_*$ is a stationary point of $f$.
\end{lemma}

\section{Line search descent method on $\St(m,n,\CC)$}\label{subsec:first_block_algo}

Let $f$ be the cost function \cref{eq:sym_tensor_diagonalization-gener-0-equi} defined on the product of $\St(m,n,\CC)$ and $\SL_{m}(\CC)$.   
Let $\matr{\omega}_{k-1}=(\matr{Y}_{k-1},\matr{X}_{k-1})$ and $p=f_{1,\matr{X}_{k-1}}$ be the first restricted function. 
Denote $\matr{X}=\matr{X}_{k-1}$ for simplicity. 
Then the restricted function $p$ can be expressed as 
\begin{equation}\label{eq:sym_tensor_diagonalization-gener-3}
p: \St(m,n,\CC)\rightarrow\RR^+,\ \ \matr{Y}\mapsto \sum_{\ell=1}^{\matrelem{L}} \|\offdiag{\matr{W}^{(\ell)}}\|^2,
\end{equation}
where 
$\matr{W}^{(\ell)}=\matr{X}^{\TH}\matr{Y}^{\TH}\matr{A}^{(\ell)}\matr{Y} \matr{X}$, and $(\cdot)^{\TH}=(\cdot)^{\T}$ or $(\cdot)^{\H}$. In this section, we adopt the \emph{line search descent} \cite{absil2005convergence,absil2009optimization,nocedal2006numerical,ring2012optimization,sato2015new} method on $\St(m,n,\CC)$ to find the next iterate $\matr{Y}_{k}$ for the restricted function $p$ in \cref{eq:sym_tensor_diagonalization-gener-3}. 

\subsection{Riemannian gradient}
We first present a lemma, which can be obtained by direct calculations.
This result will help us to obtain the Riemannian gradient of the restricted function $p$ in \cref{eq:sym_tensor_diagonalization-gener-3}. 

\begin{lemma}\label{lemm:euclid_grad_gener}
Let $\matr{A}\in\CC^{n\times n}$ and the function $\tilde{p}$ be defined as
\begin{equation*}
\tilde{p}: \CC^{n\times m}\rightarrow\RR^{+},\ \ \matr{Z}\mapsto  \|\offdiag{\matr{W}}\|^2,
\end{equation*}
where $\matr{W}=\matr{Z}^{\TH}\matr{A}\matr{Z}$. 
Denote $\matr{V}=\matr{A}\matr{Z}=[\vect{v}_1,\cdots,\vect{v}_m]\in\CC^{n\times m}$ and $\overbar{\matr{V}}=\matr{A}^{\TH}\matr{Z}=[\bar{\vect{v}}_1,\cdots,\bar{\vect{v}}_m]\in\CC^{n\times m}$. 
Denote $\matr{Z}=[\vect{z}_1,\cdots,\vect{z}_m]$. 
Then the Euclidean gradient is
\begin{align*}
\nabla\tilde{p}(\matr{Z})  = 2\left(\sum_{j\neq 1}\vect{v}_{j}\vect{v}_{j}^{\H}\vect{z}_1,\cdots,\sum_{j\neq m}\vect{v}_{j}\vect{v}_{j}^{\H}\vect{z}_m\right)+
2\left(\sum_{j\neq 1}\bar{\vect{v}}_{j}(\bar{\vect{v}}_{j})^{\H}\vect{z}_1,\cdots,\sum_{j\neq m}\bar{\vect{v}}_{j}(\bar{\vect{v}}_{j})^{\H}\vect{z}_m\right).
\end{align*}
In particular, it satisfies 
\begin{align}\label{eq:eucli_grad_h}
\matr{Z}\matr{Z}^{\H}\nabla\tilde{p}(\matr{Z})   = 2\matr{Z}\Upsilon(\matr{W}),
\end{align}
where $\Upsilon(\matr{W})\in\CC^{m\times m}$ is defined as 
\begin{equation}\label{eq:Upsilon_W}
\Upsilon(\matr{W})\eqdef\left\{
\begin{aligned}
\matr{W}\offdiag{\matr{W}}^{\H} + \matr{W}^{H}\offdiag{\matr{W}},\ \ & \textrm{if}\  (\cdot)^{\TH}=(\cdot)^{\H}; \\
\matr{W}^{*}\offdiag{\matr{W}}^{\T} + \matr{W}^{H}\offdiag{\matr{W}},\ \ & \textrm{if}\  (\cdot)^{\TH}=(\cdot)^{\T}.
\end{aligned}
\right.
\end{equation}
\end{lemma}

\begin{lemma}\label{eq:encli_grad_p}
Let $\matr{W}^{(\ell)}$ and the function $p$ be as in  \cref{eq:sym_tensor_diagonalization-gener-3}. Then
the Euclidean gradient satisfies
\begin{align}
\matr{Y}^{\H}\nabla\textit{p}(\matr{Y}) =  2(\matr{X}^{\H})^{-1}\sum_{\ell=1}^{\matrelem{L}}\Upsilon(\matr{W}^{(\ell)})\matr{X}^{\H},\label{eq:euclidean_grad-1}
\end{align}
where $\Upsilon(\matr{W}^{(\ell)})$ is as in \cref{eq:Upsilon_W}. 
\end{lemma}
\begin{proof}
By the product rule, we see that
$\nabla\textit{p}(\matr{Y})=\nabla\tilde{p}(\matr{Y}\matr{X})\matr{X}^{\H}$. 
Then, by equation \cref{eq:eucli_grad_h}, we have that 
\begin{align*}
\matr{X}\matr{X}^{\H}\matr{Y}^{\H}\nabla\textit{p}(\matr{Y}) = \matr{Y}^{\H}\matr{Y}\matr{X}(\matr{Y}\matr{X})^{\H}\nabla\tilde{p}(\matr{Y}\matr{X})\matr{X}^{\H}=
2\matr{Y}^{\H}\matr{Y}\matr{X}\sum_{\ell=1}^{\matrelem{L}}\Upsilon(\matr{W}^{(\ell)})\matr{X}^{\H}.
\end{align*}
Note that $\matr{X}$ is invertible. The proof is complete. 
\end{proof}

Now, by equations \cref{eq:euclidean_grad-1} and \cref{eq:Rie_grad}, we see that the Riemannian gradient of the function $p$ in \cref{eq:sym_tensor_diagonalization-gener-3} satisfies 
\begin{align}\label{eq:riem_grad_p}
\matr{Y}^{\H}\ProjGrad{p}{\matr{Y}} 
= (\matr{X}^{\H})^{-1}\sum_{\ell=1}^{\matrelem{L}}\Upsilon(\matr{W}^{(\ell)})\matr{X}^{\H} - \matr{X}\sum_{\ell=1}^{\matrelem{L}}\Upsilon(\matr{W}^{(\ell)})^{\H}\matr{X}^{-1}.
\end{align}

\subsection{Line search descent method}\label{subsec:line_seach}
We now begin to present more details about the line search descent method \cite{absil2005convergence,absil2009optimization,nocedal2006numerical,ring2012optimization,sato2015new} on $\St(m,n,\CC)$.  In this method, we choose the next iteration as 
\begin{equation}\label{eq:U_k_update}
\matr{Y}_{k} ={\rm Exp}_{\matr{Y}_{k-1}}(t_{k-1}\matr{V}_{k-1}),
\end{equation}
where $\matr{V}_{k-1}$ is the \emph{search direction}, $t_{k-1}$ is the \emph{step size} and ${\rm Exp}_{\matr{Y}_{k-1}}$ is the exponential map defined in \cref{eq:ex_sl_pr}. We always choose the search direction  $\matr{V}_{k-1}$ such that 
\begin{equation}\label{eq:cond_driect}
\langle\ProjGrad{p}{\matr{Y}_{k-1}},\matr{V}_{k-1}\rangle_{\matr{Y}_{k-1}}\leq -\delta_{s} \|\ProjGrad{p}{\matr{Y}_{k-1}}\|\|\matr{V}_{k-1}\|,
\end{equation}
where $0<\delta_{s}<1$ is a fixed positive constant. 
 We say that the step size $t_{k-1}$ satisfies the \emph{Armijo condition}\footnote{It is also known as the \emph{first
Wolfe condition} in the literature.}, if 
\begin{align}\label{eq:cond_step_size}
p(\matr{Y}_{k})\leq p(\matr{Y}_{k-1}) + \delta_{w}t_{k-1}\langle\ProjGrad{p}{\matr{Y}_{k-1}},\matr{V}_{k-1}\rangle_{\matr{Y}_{k-1}},
\end{align}
where $0<\delta_{w}<1$ is a fixed positive constant. We say that the step size $t_{k-1}$ satisfies the \emph{curvature condition}, if \begin{align}\label{eq:cond_step_size-2}{\small 
\langle\ProjGrad{p}{\matr{Y}_{k}},\mathbf{D}{\rm Exp}_{\matr{Y}_{k-1}}(t_{k-1}\matr{V}_{k-1})[\matr{V}_{k-1}]\rangle_{\matr{Y}_{k}}\geq \delta_{c} \langle\ProjGrad{p}{\matr{Y}_{k-1}},\matr{V}_{k-1}\rangle_{\matr{Y}_{k-1}},}
\end{align}
where $\delta_{w}<\delta_{c}<1$ is a fixed positive constant. 
The conditions \cref{eq:cond_step_size} and \cref{eq:cond_step_size-2} are known collectively as the \emph{Wolfe conditions}. 
As in the Euclidean space case \cite[Lemma 3.1]{nocedal2006numerical}, it was shown \cite{ring2012optimization,sato2015new} that we can always choose the step size $t_{k-1}$ such that the conditions \cref{eq:cond_step_size} and \cref{eq:cond_step_size-2} are both satisfied. 
It is not difficult to see that there exists $\mathrm{M}_{e}>0$ such that
\begin{equation*}
\|{\rm Exp}_{\matr{Y}}(\matr{V}_{1})-{\rm Exp}_{\matr{Y}}(\matr{V}_{2})\|\leq \mathrm{M}_{e}\|\matr{V}_{1}-\matr{V}_{2}\|,
\end{equation*}
for any $\matr{Y}\in \St(m,n,\CC)$ and $\matr{V}_1,\matr{V}_2\in{\rm\bf T}_{\matr{Y}} \St(m,n,\CC)$. 
Then the next result follows directly.
\begin{lemma}\label{lemma:line_search}
If we choose the next iterate $\matr{Y}_k$ as in \cref{eq:U_k_update} such that the conditions \cref{eq:cond_driect} and \cref{eq:cond_step_size} are both satisfied,
then we have 
\begin{align*}{\small 
p(\matr{Y}_{k-1})-p(\matr{Y}_{k})\geq \delta_{s} \delta_{w}\|\ProjGrad{p}{\matr{Y}_{k-1}}\| \|t_{k-1}\matr{V}_{k-1}\|
\geq \sigma_{p}\|\ProjGrad{p}{\matr{Y}_{k-1}}\| \|\matr{Y}_{k}-\matr{Y}_{k-1}\|,}
\end{align*}
where $\sigma_p=(\delta_{s} \delta_{w})/\mathrm{M}_{e}$. 
\end{lemma}

We also have the next result, a simple corollary of the proof in \cite[Theorem 2]{ring2012optimization}.

\begin{lemma}\label{lemma:line_search-2}
If we choose the next iterate $\matr{Y}_k$ as in \cref{eq:U_k_update} such that the conditions \cref{eq:cond_driect},  \cref{eq:cond_step_size} and \cref{eq:cond_step_size-2} are all satisfied, then we have 
\begin{equation}\label{eq:cond_weak_h}
p(\matr{Y}_{k-1})-p(\matr{Y}_{k})\geq \eta_p\|\ProjGrad{p}{\matr{Y}_{k-1}}\|^2, 
\end{equation}
where $\eta_p>0$ is a fixed positive constant. 
\end{lemma}

\section{Elementary functions and three subalgorithms}\label{sec:Jacobi_g_rota}
In this section, we define four kinds of elementary functions and present the details of three subalgorithms.

\subsection{Elementary functions and their derivatives}\label{subsec:part_deri_elem_func}

Let $g: \SL_{m}(\CC)\rightarrow\RR^+$ be a differentiable function,
 $\matr{X}\in\SL_{m}(\CC)$ and $z=x+\ui y$. 
Corresponding to the  four elementary transformations, 
we  define the following four  \emph{elementary functions}:
\begin{equation}\label{eq:elemn_func_h}
\begin{aligned}
h^{(U)}_{(i,j),\matr{X}}(x,y)&=h^{(U)}_{(i,j),\matr{X}}(\matr{\Utwo}) \eqdef g(\matr{X}\Umat{i}{j}{\matr{\Utwo}}),\ \  \matr{\Utwo} \in \SUT_{2}(\CC),\ z\in\CC;\\
h^{(L)}_{(i,j),\matr{X}}(x,y)&=h^{(L)}_{(i,j),\matr{X}}(\matr{\Utwo})  \eqdef g(\matr{X}\Lmat{i}{j}{\matr{\Utwo}}),\ \ \matr{\Utwo} \in \SLT_{2}(\CC), \ z\in\CC;\\
h^{(D)}_{(i,j),\matr{X}}(x,y)&=h^{(D)}_{(i,j),\matr{X}}(\matr{\Utwo}) \eqdef g(\matr{X}\Dmat{i}{j}{\matr{\Utwo}}),\ \ \matr{\Utwo} \in \DT_{2}(\CC), \ z\in\CC_{*};\\
h^{(Q)}_{(i,j),\matr{X}}(c,s_1,s_2)&=h^{(Q)}_{(i,j),\matr{X}}(\theta,\phi)=h^{(Q)}_{(i,j),\matr{X}}(\matr{\Utwo}) \\
&\eqdef g(\matr{X}\Gmat{i}{j}{\matr{\Utwo}}), \matr{\Utwo} \in \SUNN{2}(\CC),\  (c,s_1,s_2)\in\mathbb{S}_2,\ (\theta,\phi)\in\RR^2.
\end{aligned}
\end{equation}
In the above last equation, as in \cite{Como10:book,ULC2019}, we parameterize $\matr{\Utwo} \in \SUNN{2}(\CC)$ as  
\begin{align*}
\matr{\Psi}=\matr{\Utwo}(c,s_1,s_2)
&= \begin{bmatrix}
c & -s \\
s^{\ast} & c
\end{bmatrix}
= \begin{bmatrix}
c & -(s_1+\ui s_2) \\
s_1-\ui s_2 & c
\end{bmatrix}= \begin{bmatrix}
\cos\theta & -\sin\theta \ue^{\ui\phi} \\
\sin\theta \ue^{-\ui\phi} & \cos\theta
\end{bmatrix},
\end{align*}
where $(c,s_1,s_2)\in\mathbb{S}_{2}$ and $(\theta,\phi)\in\RR^2$.

Recall that $\matr{\Lambda}(\matr{X})\eqdef\matr{X}^{-1}\ProjGrad{g}{\matr{X}}\in\sll_{m}(\CC)$ in \Cref{subsec:grad_SL_n}, and
denote $\matr{\Lambda}=\matr{\Lambda}(\matr{X})$ for simplicity. 
We now show the relationships between the Riemannian gradients of the four elementary functions defined in \cref{eq:elemn_func_h} and the Riemannian gradient of the function $g$ at $\matr{X}\in\SL_{m}(\CC)$.
The proof is postponed to \cref{sec:ml_proofs_wellness_Jaco_part_01}. 

\begin{lemma}\label{lem:ProjGradSubmatrix}
The Riemannian gradients of the elementary functions defined in \cref{eq:elemn_func_h} at the identity matrix $\matr{I}_2$ can be expressed as follows:
\begin{align*}
(i)\ &\ProjGrad{\hij{i}{j}{\matr{X}}}{\matr{I}_2} = {\footnotesize \begin{bmatrix}
	\frac{\ui}{2}\Im\left(\matr{\Lambda}_{ii}-\matr{\Lambda}_{jj}\right) &
	\frac{1}{2}\Re\left(\matr{\Lambda}_{ij}-\matr{\Lambda}_{ji}\right)+\frac{\ui}{2}\Im\left(\matr{\Lambda}_{ij}+\matr{\Lambda}_{ji}\right)  \\
	-\frac{1}{2}\Re\left(\matr{\Lambda}_{ij}-\matr{\Lambda}_{ji}\right)+\frac{\ui}{2}\Im\left(\matr{\Lambda}_{ij}+\matr{\Lambda}_{ji}\right) &  -\frac{\ui}{2}\Im\left(\matr{\Lambda}_{ii}-\matr{\Lambda}_{jj}\right)
\end{bmatrix};}\label{eq-gradient-h}\\
(ii)\ &\ProjGrad{h^{(U)}_{(i,j),\matr{X}}}{\matr{I}_2} 
= \begin{bmatrix}
0 &
\matr{\Lambda}_{ij} \\
0  &  0
\end{bmatrix};\quad
(iii)\ \ProjGrad{h^{(L)}_{(i,j),\matr{X}}}{\matr{I}_2} 
= \begin{bmatrix}
0 &
0 \\
\matr{\Lambda}_{ji}  &  0
\end{bmatrix};\\
(iv)\ &\ProjGrad{h^{(D)}_{(i,j),\matr{X}}}{\matr{I}_2} =  
\begin{bmatrix}
\Re\left(\matr{\Lambda}_{ii}-\matr{\Lambda}_{jj}\right) &
0 \\
0  &   \Im\left(\matr{\Lambda}_{ii}-\matr{\Lambda}_{jj}\right)
\end{bmatrix}. 
\end{align*}
\end{lemma}

The following lemma can be easily obtained from \Cref{lem:ProjGradSubmatrix}.

\begin{lemma}\label{lem:ProjGradSubmatrix_2}
The partial derivatives of the elementary functions defined in \cref{eq:elemn_func_h} satisfy
\begin{align*}
(i)\ &\partial\hij{i}{j}{\matr{X}}(\matr{I}_2) \eqdef \partial\hij{i}{j}{\matr{X}}(1,0,0)
= [0,\ -\Re\left(\matr{\Lambda}_{ij} - \matr{\Lambda}_{ji}\right), \ -\Im\left(\matr{\Lambda}_{ij} + \matr{\Lambda}_{ji}\right)]^{\T};\\
(ii)\  &\partial h^{(U)}_{(i,j),\matr{X}}(\matr{I}_2) \eqdef \partial h^{(U)}_{(i,j),\matr{X}}(0,0)
= [\Re \left(\matr{\Lambda}_{ij}\right), \  \Im \left(\matr{\Lambda}_{ij}\right)]^{\T};\\
(iii)\  &\partial h^{(L)}_{(i,j),\matr{X}}(\matr{I}_2) \eqdef
\partial h^{(L)}_{(i,j),\matr{X}}(0,0)
= [\Re \left(\matr{\Lambda}_{ji}\right),  \ \Im \left(\matr{\Lambda}_{ji}\right)]^{\T};\\
(iv)\  &\partial h^{(D)}_{(i,j),\matr{X}}(\matr{I}_2) \eqdef
\partial h^{(D)}_{(i,j),\matr{X}}(0,0)
= [\Re \left(\matr{\Lambda}_{ii} -  \matr{\Lambda}_{jj}\right),  \ \Im \left(\matr{\Lambda}_{ii} -  \matr{\Lambda}_{jj}\right)]^{\T}.
\end{align*}
\end{lemma}

\subsection{Three subalgorithms}\label{subsec:Jaco_G}
%In \Cref{alg:BCD-G-RL}, 
Let $f$ be the cost function \cref{eq:sym_tensor_diagonalization-gener-0-equi}. 
Let $\matr{\omega}_{k-1}=(\matr{Y}_{k-1},\matr{X}_{k-1})$ be the $(k-1)$-th iterate produced by \Cref{alg:BCD-G-RL}, and 
$g: \SL_{m}(\CC)\rightarrow\RR$ be the restricted function $f_{2,\matr{Y}_{k-1}}$ defined as in \cref{eq:restrc_functions}. 
%In the case of \Cref{alg:jacobi-LU-G_c-1}, we take $g$ to be the function $f$ itself. 
Let $(i_k,j_k)$ be a pair of indices satisfying $1 \le i_k < j_k \le m$. 
For simplicity, we denote 
\begin{equation}\label{eq:restr_functions}
\begin{aligned}
h^{(Q)}_k &= h^{(Q)}_{(i_k,j_k),\matr{X}_{k-1}}, \quad h^{(U)}_k = h^{(U)}_{(i_k,j_k),\matr{X}_{k-1}},  \\
h^{(L)}_k &= h^{(L)}_{(i_k,j_k),\matr{X}_{k-1}}, \quad h^{(D)}_k = h^{(D)}_{(i_k,j_k),\matr{X}_{k-1}}.
\end{aligned}
\end{equation} 
Based on the three classes of elementary transformations, the subalgorithms to update $\matr{X}_k$ in \Cref{alg:BCD-G-RL} are summarized in \autoref{alg:jacobi-LU-G_c},  \autoref{alg:jacobi-LQ-G_c} and \autoref{alg:jacob-GU-G_c}, respectively. 
In these three cases, as in \Cref{subsecc:bcd-g}, 
we call \Cref{alg:BCD-G-RL} the $\emph{BCD-GLU}$, $\emph{BCD-GQU}$ and $\emph{BCD-GU}$ algorithms, respectively. 

\begin{algorithm}
\SetAlgoRefName{1a}
\SetAlgorithmName{Subalgorithm}{b}{c}
\caption{The subalgorithm to update $\matr{X}_k$ based on GLU class}\label{alg:jacobi-LU-G_c}
\begin{algorithmic}[1]
\STATE{{\bf Input:} Current iterate $\matr{X}_{k-1}$, a fixed positive constant $0<\varepsilon<\sqrt{\frac{2}{3m(m-1)}}$.}
\STATE{{\bf Output:} New iterate $\matr{X}_{k}$.}
\STATE\quad Choose an index pair $(i_k,j_k)$ and an elementary function $h_k$ such that\footnote{The inequality \cref{eq-inequality-G-LU_c} can be seen as a non-orthogonal analogue of \cite[Eq. (3.3)]{IshtAV13:simax} and \cite[Eq. (10)]{LUC2017globally}.} \begin{equation}\label{eq-inequality-G-LU_c}
\|\partial h_k(\matr{I}_2)\| 
\geq\varepsilon\|\matr{\Lambda}(\matr{X}_{k-1})\|,
\end{equation}
\quad where $h_k=h^{(U)}_k, h^{(L)}_k$ or $h^{(D)}_k$;
\STATE\quad Compute $\matr{\Utwo}^{*}_{k}$ that minimizes the elementary function $h_k$,  satisfying \Cref{remar:rho_setting_four_case_ul} and \ref{remar:rho_setting_four_case_c} (will be shown in \Cref{sec:trian_diago_rotation}); 
\STATE\quad Update $\matr{X}_k = \matr{X}_{k-1}\matr{P}_{k}$, where $\matr{P}_{k}=\Umat{i_k}{j_k}{\matr{\Utwo}^{*}_{k}}, \Lmat{i_k}{j_k}{\matr{\Utwo}^{*}_{k}}$ or $\Dmat{i_k}{j_k}{\matr{\Utwo}^{*}_{k}}$.
\end{algorithmic}
\end{algorithm}

\begin{algorithm}
\SetAlgoRefName{1b}
\SetAlgorithmName{Subalgorithm}{b}{c}
\caption{The subalgorithm to update $\matr{X}_k$ based on GQU class}\label{alg:jacobi-LQ-G_c}
\begin{algorithmic}[1]
\STATE{{\bf Input:} Current iterate $\matr{X}_{k-1}$, a fixed positive constant $0<\varepsilon<\sqrt{\frac{3-\sqrt{5}}{3m(m-1)}}$.}
\STATE{{\bf Output:} New iterate $\matr{X}_{k}$.}
\STATE\quad Choose an index pair $(i_k,j_k)$ and an elementary function $h_k$ satisfying the inequality \cref{eq-inequality-G-LU_c}, where $h_k =h^{(Q)}_k, h^{(U)}_k$ or $h^{(D)}_k$;
\STATE\quad Compute $\matr{\Utwo}^{*}_{k}$ that minimizes the elementary function $h_k$,  satisfying \Cref{remar:rho_setting_four_case_ul}, \ref{remar:rho_setting_four_case_c} and \ref{remark:condition_inner} (will be shown in \Cref{sec:trian_diago_rotation,sec:trian_diago_rotation-2}); 
\STATE\quad Update $\matr{X}_k = \matr{X}_{k-1}\matr{P}_{k}$, where $\matr{P}_{k}=\Gmat{i_k}{j_k}{\matr{\Utwo}^{*}_{k}}, \Umat{i_k}{j_k}{\matr{\Utwo}^{*}_{k}}$ or $\Dmat{i_k}{j_k}{\matr{\Utwo}^{*}_{k}}$.
\end{algorithmic}
\end{algorithm}

\begin{algorithm}
\SetAlgoRefName{1c}
\SetAlgorithmName{Subalgorithm}{b}{c}
\caption{The subalgorithm to update $\matr{X}_k$ based on GU class}\label{alg:jacob-GU-G_c}
\begin{algorithmic}[1]
\STATE{{\bf Input:} Current iterate $\matr{X}_{k-1}$, a fixed positive constant $0<\varepsilon<\sqrt{\frac{1}{m(m-1)}}$.}
\STATE{{\bf Output:} New iterate $\matr{X}_{k}$.}
\STATE\quad Choose an index pair $(i_k,j_k)$ and an elementary function $h_k$ satisfying the inequality \cref{eq-inequality-G-LU_c}, where $h_k=h^{(U)}_k$ or $h^{(D)}_k$;
\STATE\quad Compute $\matr{\Utwo}^{*}_{k}$ that minimizes the elementary function $h_k$, satisfying \Cref{remar:rho_setting_four_case_ul} and \ref{remar:rho_setting_four_case_c} (will be shown in \Cref{sec:trian_diago_rotation}); 
\STATE\quad Update $\matr{X}_k = \matr{X}_{k-1}\matr{P}_{k}$, where $\matr{P}_{k}=\Umat{i_k}{j_k}{\matr{\Utwo}^{*}_{k}}$ or $\Dmat{i_k}{j_k}{\matr{\Utwo}^{*}_{k}}$.
\end{algorithmic}
\end{algorithm}

In the following result, we will show that \autoref{alg:jacobi-LU-G_c} and \autoref{alg:jacobi-LQ-G_c} are both well-defined. The proof is postponed to \cref{sec:ml_proofs_wellness_Jaco_part_01}. 

\begin{proposition}\label{theorem-wellness-Jacobi-G-LU_c}
(i) In \autoref{alg:jacobi-LU-G_c}, 
we can always choose an index pair $(i_k,j_k)$ and an elementary function $h_k=h^{(U)}_k, h^{(L)}_k$ or $h^{(D)}_k$ such that the inequality \cref{eq-inequality-G-LU_c} is satisfied. \\ 
(ii) In \autoref{alg:jacobi-LQ-G_c}, 
we can always choose an index pair $(i_k,j_k)$ and an elementary function $h_k =h^{(Q)}_k, h^{(U)}_k$ or $h^{(D)}_k$ such that the inequality \cref{eq-inequality-G-LU_c} is satisfied. 
\end{proposition}

In BCD-GU algorithm, we always choose a starting point $\matr{X}_{0}\in\EUT_{m}(\CC)$. 
Let $\eut_{m}(\CC)\subseteq\CC^{m\times m}$ be the set of upper triangular matrices with the trace equal to 0. 
Then the tangent space to $\EUT_{m}(\CC)$ at a point $\matr{X}\in\EUT_{m}(\CC)$ can be constructed \cite{afsari2004gradient,baker2012matrix} by $\mathbf{T}_{\matr{X}}\EUT_{m}(\CC) = \{\matr{X}\Omega, \Omega\in\eut_{m}(\CC)\}$, which is useful to the proof of the following result. The proof is postponed to \cref{sec:ml_proofs_wellness_Jaco_part_01}. 

\begin{proposition}\label{theorem-wellness-Jacobi-G-GU_c}
In \autoref{alg:jacob-GU-G_c},
we can always choose an index pair $(i_k,j_k)$ and an elementary function $h_k=h^{(U)}_k$ or $h^{(D)}_k$ such that the inequality \cref{eq-inequality-G-LU_c} is satisfied. 
\end{proposition}

\section{Plane triangular and diagonal transformations for JADM problem}\label{sec:trian_diago_rotation}

Let $f$ be the cost function \cref{eq:sym_tensor_diagonalization-gener-0-equi}. 
Let $\matr{\omega}_{k-1}=(\matr{Y}_{k-1},\matr{X}_{k-1})$ and $g: \SL_{m}(\CC)\rightarrow\RR$ be the restricted function $f_{2,\matr{Y}_{k-1}}$ as in \Cref{subsec:Jaco_G}. 
Denote $\matr{B}^{(\ell)}=\matr{Y}_{k-1}^{\TH}\matr{A}^{(\ell)} \matr{Y}_{k-1}$ for $1\leq\ell\leq \matrelem{L}$, where $(\cdot)^{\TH}=(\cdot)^{\T}$ or $(\cdot)^{\H}$. 
Then $g$ can be expressed as 
\begin{equation}\label{eq:sym_tensor_diagonalization-gener-4}
g: \SL_{m}(\CC)\rightarrow\RR^+,\ \ \matr{X}\mapsto \sum_{\ell=1}^{\matrelem{L}} \|\offdiag{\matr{W}^{(\ell)}}\|^2,
\end{equation}
where 
$\matr{W}^{(\ell)}=\matr{X}^{\TH}\matr{B}^{(\ell)} \matr{X}$ for $1\leq\ell\leq \matrelem{L}$. 
In this section, we will first calculate the Riemannian gradient of $g$ in \cref{eq:sym_tensor_diagonalization-gener-4}, and the partial derivatives of elementary functions $h^{(U)}_k$, $h^{(L)}_k$ and $h^{(D)}_k$ in \cref{eq:restr_functions}. 
Then, we will prove that inequalities \cref{eq:sufficient_descent-2} and \cref{equation-condition-weak_c} are both satisfied in the plane triangular and diagonal transformations.

\subsection{Riemannian gradient}

Let $g$ and $\matr{W}^{(\ell)}$ be as in  \cref{eq:sym_tensor_diagonalization-gener-4}.  
Then, by equations \cref{eq:eucli_grad_h} and \cref{eq:riema_gradient_c}, 
we have the Euclidean gradient and Riemannian gradient of $g$ at $\matr{X}\in\SL_{m}(\CC)$ as follows: 
\begin{align}
\nabla\textit{g}(\matr{X})  & = 2(\matr{X}^{\H})^{-1}\sum_{\ell=1}^{\matrelem{L}}\Upsilon(\matr{W}^{(\ell)}),\label{eq:euclidean_grad}\\
\ProjGrad{g}{\matr{X}} & = 2\matr{X}\sum_{\ell=1}^{\matrelem{L}}\left(\Upsilon(\matr{W}^{(\ell)}) - \frac{\tr(\Upsilon(\matr{W}^{(\ell)}))}{n}\matr{I}_{n}\right),\label{eq:rrieman_grad}
\end{align}
where $\Upsilon(\matr{W}^{(\ell)})$ is defined as in equation \cref{eq:Upsilon_W}. 

\begin{remark}
In the real case, the Euclidean gradient in \cref{eq:euclidean_grad} was earlier derived in \cite[Eq. (6.3)]{afsari2004gradient} and \cite[Section 2.3]{bouchard2020approximate}.
In this paper, we extend it to problem \cref{eq:sym_tensor_diagonalization-gener-4} in the complex case and calculate the Riemannian gradient \cref{eq:rrieman_grad} as well. 
\end{remark}

\subsection{Elementary functions}\label{subsec:ele_func_trian}

Let $\matr{W}^{(\ell)}=\matr{X}_{k-1}^{\TH}\matr{B}^{(\ell)} \matr{X}_{k-1}$ for $1\leq\ell\leq \matrelem{L}$. 
Let
\begin{equation}\label{eq:varrho_1}
\varrho\eqdef\left\{
\begin{aligned}
1,\ \ & \textrm{if}\  (\cdot)^{\TH}=(\cdot)^{\H}; \\
-1,\ \ & \textrm{if}\  (\cdot)^{\TH}=(\cdot)^{\T}.
\end{aligned}
\right.
\end{equation}
Denote $(i,j)=(i_k,j_k)$ for simplicity. 
Now we use the following notations:
\begin{itemize}
\item 
$\begin{aligned}[t]
&\alpha_{1} \eqdef  \sum_{\ell=1}^{\matrelem{L}}\sum_{p\neq j}\left(|\matrelem{W}^{(\ell)}_{ip}|^2+|\matrelem{W}^{(\ell)}_{pi}|^2\right), \\
&\alpha_{2} \eqdef  \sum_{\ell=1}^{\matrelem{L}}\sum_{p\neq j}\left(\matrelem{W}^{(\ell,\R)}_{ip}\matrelem{W}^{(\ell,\R)}_{jp}+\matrelem{W}^{(\ell,\I)}_{ip}\matrelem{W}^{(\ell,\I)}_{jp}+\matrelem{W}^{(\ell,\R)}_{pi}\matrelem{W}^{(\ell,\R)}_{pj}+\matrelem{W}^{(\ell,\I)}_{pi}\matrelem{W}^{(\ell,\I)}_{pj}\right), \\
&\alpha_{3} \eqdef   \sum_{\ell=1}^{\matrelem{L}}\sum_{p\neq j}\left(\varrho\left(\matrelem{W}^{(\ell,\I)}_{ip}\matrelem{W}^{(\ell,\R)}_{jp}-\matrelem{W}^{(\ell,\R)}_{ip}\matrelem{W}^{(\ell,\I)}_{jp}\right)+\matrelem{W}^{(\ell,\R)}_{pi}\matrelem{W}^{(\ell,\I)}_{pj}-\matrelem{W}^{(\ell,\I)}_{pi}\matrelem{W}^{(\ell,\R)}_{pj}\right).
\end{aligned}$
\item 
$\begin{aligned}[t]
&\beta_{1} \eqdef  \sum_{\ell=1}^{\matrelem{L}}\sum_{p\neq i}\left(|\matrelem{W}^{(\ell)}_{jp}|^2+|\matrelem{W}^{(\ell)}_{pj}|^2\right), \\
&\beta_{2} \eqdef  \sum_{\ell=1}^{\matrelem{L}}\sum_{p\neq i}\left(\matrelem{W}^{(\ell,\R)}_{ip}\matrelem{W}^{(\ell,\R)}_{jp}+\matrelem{W}^{(\ell,\I)}_{ip}\matrelem{W}^{(\ell,\I)}_{jp}+\matrelem{W}^{(\ell,\R)}_{pi}\matrelem{W}^{(\ell,\R)}_{pj}+\matrelem{W}^{(\ell,\I)}_{pi}\matrelem{W}^{(\ell,\I)}_{pj}\right),  \\
&\beta_{3} \eqdef  \sum_{\ell=1}^{\matrelem{L}}\sum_{p\neq i}\left(\varrho\left(\matrelem{W}^{(\ell,\R)}_{ip}\matrelem{W}^{(\ell,\I)}_{jp}-\matrelem{W}^{(\ell,\I)}_{ip}\matrelem{W}^{(\ell,\R)}_{jp}\right)+\matrelem{W}^{(\ell,\I)}_{pi}\matrelem{W}^{(\ell,\R)}_{pj}-\matrelem{W}^{(\ell,\R)}_{pi}\matrelem{W}^{(\ell,\I)}_{pj}\right).
\end{aligned}$
\item 
$\begin{aligned}[t]
&\gamma_{1} \eqdef  \sum_{\ell=1}^{\matrelem{L}}\sum_{p\neq i,j}\left(|{\matrelem{W}^{(\ell)}_{ip}}|^2+|{\matrelem{W}^{(\ell)}_{pi}}|^2\right), \ \ \ 
\gamma_{2} \eqdef  \sum_{\ell=1}^{\matrelem{L}}\sum_{p\neq i,j}\left(|{\matrelem{W}^{(\ell)}_{jp}}|^2+|{\matrelem{W}^{(\ell)}_{pj}}|^2\right). 
\end{aligned}$
\end{itemize}
Then we can get the following results by direct calculations.

\begin{lemma}\label{lemm:elemen_express_c}
Let the function $g$ be as in \cref{eq:sym_tensor_diagonalization-gener-4}. Then\\
(i) the elementary function $h^{(U)}_{k}$ in \cref{eq:restr_functions} and its optimal solution $(x_k^*,y_k^*)$ satisfy
\begin{align}
h^{(U)}_{k}(x,y) - h^{(U)}_{k}(0,0) & = \alpha_{1}x^2 + 2\alpha_{2}x + \alpha_{1}y^2 + 2 \alpha_{3}y, \notag\\
h^{(U)}_{k}(x_k^*,y_k^*)-h^{(U)}_{k}(0,0)  & =  -\frac{1}{\alpha_{1}}\left(\alpha_{2}^2+\alpha_{3}^2\right), \label{eq:elemen_matr_01_c_s}\\
\partial h^{(U)}_{k}(0,0) 
& = 2[\alpha_{2},\  \alpha_{3}]^{\T}. \notag
\end{align}
(ii) the elementary function $h^{(L)}_k$ in \cref{eq:restr_functions} and its optimal solution $(x_k^*,y_k^*)$ satisfy 
\begin{align}
h^{(L)}_{k}(x,y) - h^{(L)}_{k}(0,0) 
& = \beta_{1}x^2 + 2\beta_{2}x + \beta_{1}y^2 + 2\beta_{3}y, \notag\\
h^{(L)}_{k}(x_k^*,y_k^*)-h^{(L)}_{k}(0,0) & = -\frac{1}{\beta_{1}}\left(\beta_{2}^2+\beta_{3}^2\right), \notag\\
\partial h^{(L)}_{k}(0,0) & = 2[\beta_{2},\ \beta_{3}]^{\T}. \notag
\end{align}
(iii) the elementary function $h^{(D)}_k$ in \cref{eq:restr_functions} and its optimal solution $(x_k^*,y_k^*)$ satisfy  
\begin{align}
h^{(D)}_{k}(x,y) - h^{(D)}_{k}(1,0) & = \gamma_{1}(x^2+y^2) + \gamma_{2}\frac{1}{x^2+y^2} -\gamma_{1} - \gamma_{2},\notag\\
h^{(D)}_{k}(x_k^*,y_k^*) - h^{(D)}_{k}(1,0)
& = -\left(\sqrt{\gamma_{1}}-\sqrt{\gamma_{2}}\right)^2,\notag\\
\partial h^{(D)}_{k}(1,0) & = 2[\gamma_{1} - \gamma_{2},\ 0]^{\T} .\notag
\end{align}
\end{lemma}

\begin{remark}
In the real case,
the solution $x_k^{*}$ in \cref{eq:elemen_matr_01_c_s} was earlier derived in \cite[Eq. (7)]{afsari2006simple}.
In the complex case,
the solution $z_k^{*}=x_k^{*}+\ui y_k^{*}$ in \cref{eq:elemen_matr_01_c_s} was earlier derived in 
\cite[Eq. (8)]{wang2012complex}.
\end{remark}

\begin{setting}\label{remar:rho_setting_four_case_ul}
In \Cref{alg:BCD-G-RL} for cost function \cref{eq:sym_tensor_diagonalization-gener-0-equi},
when the elementary function $h_k=h^{(U)}_{k}$,
we see that $x_{k}^{*}=0$ if $\alpha_{1}\neq0$ and $\alpha_{2}=0$. 
It is not possible that $\alpha_{1}=0$ and $\alpha_{2}\neq0$. 
If $\alpha_{1}=\alpha_{2}=0$, we set $x_{k}^{*}=0$. 
In the case of $h_k=h^{(L)}_{k}$,
we make the similar update rules for the value of $y_{k}^{*}$. 
\end{setting}

\begin{setting}\label{remar:rho_setting_four_case_c}
Let $0<\varsigma_{D}<\frac{1}{4}$ be a small positive constant.
In \Cref{alg:BCD-G-RL} for cost function \cref{eq:sym_tensor_diagonalization-gener-0-equi},
if $h_k=h^{(D)}_{k}$, we always set $y_k^{*}=0$.
Moreover, we determine $x_k^{*}$ based on the following rules.
\begin{itemize}
\item 
If $\gamma_{1}=\gamma_{2} = 0$, we set $x_k^{*}=0$.\\
\item 
Let $\varpi\eqdef\frac{\gamma_{2}}{\gamma_{1}}$. If $\varpi\in[0,\varsigma_{D})$, we set $x_k^{*}=\frac{1}{2}$.
If $\varpi\in(\frac{1}{\varsigma_{D}},+\infty]$, we set $x_k^{*}=2$.\\
\item 
Otherwise, if $\varpi\in[\varsigma_{D},\frac{1}{\varsigma_{D}}]$, we set $x_k^{*}=\sqrt[\leftroot{-2}\uproot{4}4]{\varpi}$,
which is the minimum point.
\end{itemize}
\end{setting}

\subsection{Inequalities for global convergence}
It will be seen that $f(\matr{\omega}_{k})\leq f(\matr{\omega}_{k-1})$ always holds in \Cref{alg:BCD-G-RL}.
We denote $\mathrm{M}_0\eqdef f(\matr{\omega}_0)$ in \Cref{alg:BCD-G-RL} for cost function \cref{eq:sym_tensor_diagonalization-gener-0-equi}. 
Then we have that $\gamma_{1}+\gamma_{2}\leq\mathrm{M}_0=f(\matr{\omega}_0)$. 
In the following result, we will show an inequality, which is helpful to establish inequality \cref{eq:sufficient_descent-2} when elementary function $h_k=h^{(D)}_{k}$. 
The proof is postponed to
\cref{sec:proofs_section_5}.

\begin{lemma}\label{lemma:global_rho_c-0}
In \Cref{alg:BCD-G-RL} for cost function \cref{eq:sym_tensor_diagonalization-gener-0-equi},
there exists $\iota_{D}>0$ such that 
\begin{equation}\label{equation-condition-weak_c-5}
g(\matr{X}_{k-1})-g(\matr{X}_{k})\geq\iota_{D}\|\matr{\Lambda}(\matr{X}_{k-1})\|\|\matr{\Utwo}^{*}_{k}-\matr{I}_{2}\|,
\end{equation}
whenever the elementary function $h_k=h^{(D)}_{k}$. 	
\end{lemma}

Note that 
\begin{align}\label{eq:x_k_x_k-1_1}
\|\matr{X}_{k}-\matr{X}_{k-1}\|&\leq\|\matr{\Utwo}^{*}_{k}-\matr{I}_{2}\|\|\matr{X}_{k-1}\|,\\ \|\ProjGrad{g}{\matr{X}_{k-1}}\|&\leq\|\matr{\Lambda}(\matr{X}_{k-1})\|\|\matr{X}_{k-1}\|.\label{eq:x_k_x_k-1_2}
\end{align}
Let $\iota_{D}$ be as in \cref{equation-condition-weak_c-5}, and $\mathrm{M}_{\omega}$ be as in the condition \cref{eq:condi_bounded_X_gene}. 
Let $\sigma_{D}=\iota_{D}/\mathrm{M}_{\omega}^2>0$. 
Then the next result follows directly from \Cref{lemma:global_rho_c-0}, inequalities \cref{eq:x_k_x_k-1_1} and \cref{eq:x_k_x_k-1_2}. 

\begin{corollary}\label{lemma:global_rho_c}
In \Cref{alg:BCD-G-RL} for cost function \cref{eq:sym_tensor_diagonalization-gener-0-equi},
if the iterates remain bounded, i.e., the condition \cref{eq:condi_bounded_X_gene} is satisfied, then 
\begin{equation}\label{equation-condition-weak_c-5-0}
g(\matr{X}_{k-1})-g(\matr{X}_{k})\geq\sigma_{D}\|\ProjGrad{g}{\matr{X}_{k-1}}\|\|\matr{X}_{k}-\matr{X}_{k-1}\|,
\end{equation}
whenever the elementary function $h_k=h^{(D)}_{k}$. 	
\end{corollary}

As for inequality \cref{equation-condition-weak_c-5}, we now show a similar result for the cases of elementary functions $h^{(L)}_{k}$ and $h^{(U)}_k$. 
The proof is also postponed to
\cref{sec:proofs_section_5}. 

\begin{lemma}\label{lemma:global_rho_c-2-0}
In \Cref{alg:BCD-G-RL} for cost function \cref{eq:sym_tensor_diagonalization-gener-0-equi},
there exists $\iota_{LU}>0$ such that 
\begin{equation}\label{equation-condition-weak_c-3}
g(\matr{X}_{k-1})-g(\matr{X}_{k})\geq\iota_{LU}\|\matr{\Lambda}(\matr{X}_{k-1})\|\|\matr{\Utwo}^{*}_{k}-\matr{I}_{2}\|,
\end{equation}
whenever the elementary function $h_k=h^{(L)}_{k}$ or $h^{(U)}_k$. 
\end{lemma}

Let $\iota_{LU}$ be as in \cref{equation-condition-weak_c-3} and 
$\sigma_{LU}=\iota_{LU}/\mathrm{M}_{\omega}^2>0$. 
Similar as for \Cref{lemma:global_rho_c}, the next result follows directly from \Cref{lemma:global_rho_c-2-0}, inequalities \cref{eq:x_k_x_k-1_1} and \cref{eq:x_k_x_k-1_2}.  

\begin{corollary}\label{lemma:global_rho_c-2}
In \Cref{alg:BCD-G-RL} for cost function \cref{eq:sym_tensor_diagonalization-gener-0-equi},
if the iterates remain bounded, i.e., the condition \cref{eq:condi_bounded_X_gene} is satisfied, then 
\begin{equation*}
g(\matr{X}_{k-1})-g(\matr{X}_{k})\geq\sigma_{LU}\|\ProjGrad{g}{\matr{X}_{k-1}}\|\|\matr{X}_{k}-\matr{X}_{k-1}\|,
\end{equation*}
whenever the elementary function $h_k=h^{(L)}_{k}$ or $h^{(U)}_k$.  	
\end{corollary}

\subsection{Inequalities for weak convergence}
In this subsection, we show an inequality, which will be helpful to establish inequality \cref{equation-condition-weak_c}. 
The proof is postponed to
\cref{sec:proofs_section_5}. 

\begin{lemma}\label{lemma:weak_rho_inequality_c-1}
In \Cref{alg:BCD-G-RL} for cost function \cref{eq:sym_tensor_diagonalization-gener-0-equi},
if the iterates remain bounded, i.e., the condition \cref{eq:condi_bounded_X_gene} is satisfied, then there exists $\kappa>0$ such that 
\begin{equation*}
\|\matr{\Utwo}^{*}_{k}-\matr{I}_{2}\|\geq	\kappa\|\matr{\Lambda}(\matr{X}_{k-1})\|,
\end{equation*}
whenever the elementary function $h_k=h^{(D)}_{k}$, $h^{(L)}_{k}$ or $h^{(U)}_k$. 
\end{lemma}

By \Cref{lemma:global_rho_c-0}, \Cref{lemma:global_rho_c-2-0} and \Cref{lemma:weak_rho_inequality_c-1}, we can easily get the following results by setting $\eta_{D}=(\kappa\iota_{D})/\mathrm{M}_{\omega}^2$ and 
$\eta_{LU}=(\kappa\iota_{LU})/\mathrm{M}_{\omega}^2$.

\begin{corollary}\label{lemma:weak_rho_inequality_c}
In \Cref{alg:BCD-G-RL} for cost function \cref{eq:sym_tensor_diagonalization-gener-0-equi},
if the iterates remain bounded, i.e., the condition \cref{eq:condi_bounded_X_gene} is satisfied, then 
\begin{equation}\label{equation-condition-weak_c-4}
g(\matr{X}_{k-1})-g(\matr{X}_{k})\geq\min(\eta_{D},\eta_{LU})\|\ProjGrad{g}{\matr{X}_{k-1}}\|^2,
\end{equation}
whenever the elementary function $h_k=h^{(D)}_{k}$, $h^{(L)}_{k}$ or $h^{(U)}_k$. 
\end{corollary}

\section{Givens plane transformations for JADM problem}\label{sec:trian_diago_rotation-2}

Let the function $g$ be as in \cref{eq:sym_tensor_diagonalization-gener-4}.
Let $\matr{W}^{(\ell)}=\matr{X}_{k-1}^{\TH}\matr{B}^{(\ell)} \matr{X}_{k-1}$ for $1\leq\ell\leq \matrelem{L}$ as in \Cref{subsec:ele_func_trian}, where $(\cdot)^{\TH}=(\cdot)^{\T}$ or $(\cdot)^{\H}$. 
Let $\varrho$ be as in \cref{eq:varrho_1}. 
Denote $(i,j)=(i_k,j_k)$ for simplicity. 
Define 
\begin{align}\label{eq:gamma_matr_ij}
\Gamij{i}{j}{\matr{X}_{k-1}}\eqdef
\frac{\varrho}{2}\sum\limits_{\ell=1}^{L}   \R\left(\vect{z}_{i,j}(\matr{W}^{(\ell)}) \vect{z}_{i,j}^{\H}(\matr{W}^{(\ell)})\right)\in\RR^{3\times 3},
\end{align}
where
\begin{align*}
\vect{z}_{i,j}(\matr{W})\eqdef\left\{
\begin{aligned}
&\begin{bmatrix} \matrelem{W}_{jj}-\matrelem{W}_{ii}, & \matrelem{W}_{ij}+\matrelem{W}_{ji}, & -\ui (\matrelem{W}_{ij}-\matrelem{W}_{ji}) \end{bmatrix}^{\T},\ \  &\textrm{if}\  (\cdot)^{\TH}=(\cdot)^{\H}; \\
&\begin{bmatrix} \matrelem{W}_{ij}+\matrelem{W}_{ji}, & \matrelem{W}_{ii}-\matrelem{W}_{jj}, & \ui (\matrelem{W}_{ii}+\matrelem{W}_{jj}) \end{bmatrix}^{\T},\ \  &\textrm{if}\  (\cdot)^{\TH}=(\cdot)^{\T}.
\end{aligned}
\right.
\end{align*}
Denote 
\begin{align*}
c_0 \eqdef 
\begin{cases} \frac{1}{2}\sum_{\ell=1}^{\matrelem{L}} \left|\matrelem{W}^{(\ell)}_{jj}-\matrelem{W}^{(\ell)}_{ii}\right|^2,\ \  &\textrm{if}\  (\cdot)^{\TH}=(\cdot)^{\H}; \\  -\frac{1}{2}\sum_{\ell=1}^{\matrelem{L}} \left|\matrelem{W}^{(\ell)}_{ij}+\matrelem{W}^{(\ell)}_{ji}\right|^2,\ \  &\textrm{if}\  (\cdot)^{\TH}=(\cdot)^{\T}.
\end{cases} 
\end{align*}

\subsection{Elementary function}\label{subsec:eleme_func_H}
As in \cite[Eq. (4.4)]{ULC2019}, we denote the unit vector 
\begin{align}
\vect{r}  
\eqdef\left[2c^2-1,\ -2cs_1,\ -2cs_2\right]^{\T} 
= \left[\cos 2\theta,\ -\sin 2\theta\cos\phi,\  -\sin 2\theta\sin\phi\right]^{\T}, \label{eq:w_definition}
\end{align}
where $c\in \RR^{+}$, $s= s_1 + \ui s_2\in\CC$, $c^2+|s|^2 =1$, and $\theta,\phi\in\RR$ are two angles. 
Then we can get the following result\footnote{In the $(\cdot)^{\TH}=(\cdot)^{\H}$ case, this expression was first formulated in \cite{cardoso1996jacobi}.} by direct calculations.

\begin{lemma}\label{lemm:elemen_express_c_1}
In \Cref{alg:BCD-G-RL} for cost function \cref{eq:sym_tensor_diagonalization-gener-0-equi},
the elementary function $h^{(Q)}_k$ satisfies 
\begin{align}\label{eq:elemen_h_Gamma}
h^{(Q)}_k(c,s_1,s_2) - h^{(Q)}_k(1,0,0) =
-\left(\vect{r}^{\T}\Gamij{i}{j}{\matr{X}_{k-1}}\vect{r} - c_0\right),
\end{align}
where $\Gamij{i}{j}{\matr{X}_{k-1}}\in\RR^{3\times 3}$ is as in equation \cref{eq:gamma_matr_ij}.
\end{lemma}

Denote $\matr{\Gamma}=\Gamij{i}{j}{\matr{X}_{k-1}}$ for simplicity. 
It follows by equations \cref{eq:w_definition} and \cref{eq:elemen_h_Gamma} that
\begin{align}
h^{(Q)}_k(c,s_1,s_2) - h^{(Q)}_k(1,0,0) &= -\left(q(\theta,\phi)-c_0\right),\label{eq:h_k_theta}
\end{align}
where {\small
\begin{align}\label{eq:q_tehta}
q(\theta,\phi)&\eqdef  \frac{1}{2}\left(\Gamma_{11}-\Gamma_{22}\cos^2\phi-\Gamma_{33}\sin^2\phi-\Gamma_{23}\sin(2\phi)\right)\cos(4\theta)\notag\\ &\quad-\left(\Gamma_{12}\cos\phi+\Gamma_{13}\sin\phi\right)\sin(4\theta)
+\frac{1}{2}\left(\Gamma_{11}+\Gamma_{22}\cos^2\phi+\Gamma_{33}\sin^2\phi+\Gamma_{23}\sin(2\phi)\right).
\end{align}}
Note that, by \Cref{lemm:elemen_express_c_1} and equation \cref{eq:w_definition}, we have
\begin{align}\label{eq:partial_h_k}
\partial h^{(Q)}_k(\matr{I}_2) 
=-4 [0,\ \Gamma_{12},\ \Gamma_{13}]^\T.
\end{align}

\begin{remark}
By equation \cref{eq:h_k_theta}, we see that $h^{(Q)}_k(\theta+\pi/{2},\phi)=h^{(Q)}_k(\theta,\phi)$ for any $\theta,\phi\in\RR$.
Therefore, we can always choose  $\theta_{*}\in[-\pi/{4},\pi/{4}]$. 
\end{remark}

\begin{setting}\label{remark:condition_inner}
In \Cref{alg:BCD-G-RL} for cost function \cref{eq:sym_tensor_diagonalization-gener-0-equi},
we set a positive constant $\varsigma_{Q}>0$. 
If $h_k=h^{(Q)}_k$,
we find the eigenvector $\vect{u}$ of $\matr{\Gamma}$ corresponding to the largest eigenvalue. 
Define two vectors  
$\vect{v}_{i,j}\eqdef\left[\Gamma_{12},\ \Gamma_{13}\right]^\T\in\RR^2$
and $\vect{w}_{i,j}\eqdef
\left[u_2,\ u_3\right]^\T\in\RR^2.$
\begin{itemize}
\item If it holds that 
\begin{equation}\label{eq:condition_phi_1}
|\langle\vect{v}_{i,j},\vect{w}_{i,j}\rangle|
\geq\varsigma_{Q}\|\vect{v}_{i,j}\|\|\vect{w}_{i,j}\|,
\end{equation}
then we find $\phi_*$ and $\theta_*$ by setting $\vect{r}=\vect{u}$, and $\matr{\Utwo}^{*}_{k}=\matr{\Utwo}(\theta_*,\phi_*)$;
\item Otherwise, we set $[\cos\phi_*,\ \sin\phi_*]^\T=\vect{v}_{i,j}/{\|\vect{v}_{i,j}\|}$, and then calculate $\theta_*$, which maximizes the restricted function $q(\theta,\phi_{*})$. 
\end{itemize}
\end{setting}

\subsection{Inequalities for global convergence}\label{subsec:condi_glob}

We first present a lemma, which will help us to prove \Cref{lem:p_theta_phi}.

\begin{lemma}\label{lemm:g-theta}
Let $\alpha,\beta\in\RR$ be two constants. 
For $\theta\in[-\frac{\pi}{4},\frac{\pi}{4}]$, we define a function $p(\theta)\eqdef \alpha\cos(4\theta) + \beta\sin(4\theta)$. 
If $\theta_{*}\in[-\frac{\pi}{4},\frac{\pi}{4}]$ satisfies $p(\theta_{*})=\max p(\theta)$, then we have 
\begin{equation*}
p(\theta_{*})-p(0)\geq 2\sqrt{2}|\beta|\left|\sin(\frac{\theta_*}{2})\right|.
\end{equation*} 
\end{lemma}

\begin{lemma}\label{lem:p_theta_phi}
Let the function $q(\theta,\phi)$ be as in equation \cref{eq:q_tehta}.
Suppose that $\phi_*$ and $\theta_*$ are determined as in \Cref{remark:condition_inner}. 
Then we have
\begin{equation*}
q(\theta_*,\phi_*)-q(0,0)\geq 2\sqrt{2}\varsigma_{Q} \left|\sin(\frac{\theta_{*}}{2})\right|\|\vect{v}_{i,j}\|, 
\end{equation*}
where $\varsigma_{Q}$ is the positive constant defined in \Cref{remark:condition_inner}. 
\end{lemma}
\begin{proof}
By \Cref{remark:condition_inner}, we see that 
\begin{equation}\label{eq:condition_phi_2}
|\langle\vect{v}_{i,j},[\cos\phi_*\ \sin\phi_*]^\T\rangle|
\geq\varsigma_{Q}\|\vect{v}_{i,j}\|
\end{equation}
always holds. 
By \Cref{lemm:g-theta} and the above inequality \cref{eq:condition_phi_2}, we get that 
\begin{align*}
q(\theta_{*},\phi_{*})-q(0,0)&=q(\theta_{*},\phi_{*})-q(0,\phi_{*})\geq 2\sqrt{2}\left|\sin(\frac{\theta_*}{2})\right||\Gamma_{12}\cos\phi_*+\Gamma_{13}\sin\phi_*|\\
&\geq 2\sqrt{2}
\varsigma_{Q}\left|\sin(\frac{\theta_*}{2})\right|\|\vect{v}_{i,j}\|.
\end{align*}
The proof is complete.	
\end{proof}

As for inequality \cref{equation-condition-weak_c-5}, we now show a similar result for the case of elementary functions $h^{(Q)}_{k}$, which will be helpful to establish inequality \cref{eq:sufficient_descent-2}. 

\begin{lemma}\label{lemma:h_k_condition_SU-0}
In \Cref{alg:BCD-G-RL} for cost function \cref{eq:sym_tensor_diagonalization-gener-0-equi},
there exists $\iota_{Q}>0$ such that \begin{equation}\label{equation-condition-weak_c-5-1}
g(\matr{X}_{k-1})-g(\matr{X}_{k})\geq\iota_{Q}\|\matr{\Lambda}(\matr{X}_{k-1})\|\|\matr{\Utwo}^{*}_{k}-\matr{I}_{2}\|,
\end{equation}
whenever the elementary function $h_k=h^{(Q)}_{k}$. 
\end{lemma}
\begin{proof}
We only prove the case $(\cdot)^{\TH}=(\cdot)^{\H}$, the other case being similar.  
By \Cref{lem:p_theta_phi} and equation \cref{eq:partial_h_k}, we get that
\begin{align*}
h^{(Q)}_k(0,0)-h^{(Q)}_k(\theta_*,\phi_*) &\geq 2\sqrt{2}\varsigma_{Q} \left|\sin(\frac{\theta_{*}}{2})\right|\|\vect{v}_{i,j}\|
= \frac{\varsigma_{Q}}{4} 2\sqrt{2}\left|\sin(\frac{\theta_{*}}{2})\right|\|\partial h^{(Q)}_k(\matr{I}_2)\|\\
&\geq \frac{\varsigma_{Q}\varepsilon}{4}  \|\Gmat{i}{j}{\matr{\Utwo}_{k}^{*}}-\matr{I}_m\|\|\matr{\Lambda}(\matr{X}_{k-1})\|.
\end{align*}
We can set $\iota_{Q}=(\varsigma_{Q}\varepsilon)/{4}$. The proof is complete. 
\end{proof}

Let $\iota_{Q}$ be as in \cref{equation-condition-weak_c-5-1} and $\mathrm{M}_{\omega}$ be as in condition  \cref{eq:condi_bounded_X_gene}. 
Let $\sigma_{Q}=\iota_{Q}/\mathrm{M}_{\omega}^2>0$. 
As for \Cref{lemma:global_rho_c}, the next result follows directly from \Cref{lemma:h_k_condition_SU-0}, inequalities \cref{eq:x_k_x_k-1_1} and \cref{eq:x_k_x_k-1_2}.

\begin{corollary}\label{lemma:h_k_condition_SU}
In \Cref{alg:BCD-G-RL} for cost function \cref{eq:sym_tensor_diagonalization-gener-0-equi},
if the iterates remain bounded, i.e., the condition \cref{eq:condi_bounded_X_gene} is satisfied,
then 
\begin{equation*}
g(\matr{X}_{k-1})-g(\matr{X}_{k})\geq\sigma_{Q}\|\ProjGrad{g}{\matr{X}_{k-1}}\|\|\matr{X}_{k}-\matr{X}_{k-1}\|,
\end{equation*}
whenever the elementary function $h_k=h^{(Q)}_{k}$. 
\end{corollary}

\subsection{Inequalities for weak convergence}
If the condition \cref{eq:condi_bounded_X_gene} is satisfied, it is easy to see that there exists a positive constant  $\mathrm{M}_{\Gamma}>0$ such that $\|\Gamij{i}{j}{\matr{X}_{k-1}}\|\leq \mathrm{M}_{\Gamma}$ always holds in \Cref{alg:BCD-G-RL}. 
In this subsection, we first show an inequality, which will be helpful to establish inequality \cref{equation-condition-weak_c}. 

\begin{lemma}\label{lemma:weak_rho_inequality_c-3}
In \Cref{alg:BCD-G-RL} for cost function \cref{eq:sym_tensor_diagonalization-gener-0-equi},
if the iterates remain bounded, i.e., the condition \cref{eq:condi_bounded_X_gene} is satisfied, then there exists $\kappa_{Q}>0$ such that 
\begin{equation}\label{equation-condition-weak_c-4-3}
\|\matr{\Utwo}^{*}_{k}-\matr{I}_{2}\|\geq	\kappa_{Q}\|\matr{\Lambda}(\matr{X}_{k-1})\|,
\end{equation}
whenever the elementary function $h_k=h^{(Q)}_{k}$. 
\end{lemma}

\begin{proof}
If $\vect{v}_{i,j}$ and $\vect{w}_{i,j}$ satisfy  inequality \cref{eq:condition_phi_1}, then  inequality \cref{equation-condition-weak_c-4-3} can be proved by a similar method as for \cite[Lemma 7.2]{ULC2019}.
Otherwise, if we set $[\cos\phi_*\ \sin\phi_*]^\T=\vect{v}_{i,j}/{\|\vect{v}_{i,j}\|}$ and find $\theta_*$ based on $\phi_*$, then 
\begin{align*}
|\sin(4\theta_{*})|&=\frac{\left|\Gamma_{12}\cos\phi_*+\Gamma_{13}\sin\phi_*\right|}{\sqrt{\left(\Gamma_{12}\cos\phi_*+\Gamma_{13}\sin\phi_*\right)^2+\frac{1}{4}\left(\Gamma_{11}-\Gamma_{22}\cos^2\phi_*-\Gamma_{33}\sin^2\phi_*-\Gamma_{23}\sin(2\phi_*)\right)^2}}\\
&\geq \frac{\sqrt{\Gamma_{12}^2+\Gamma_{13}^2}}{2\sqrt{5}\mathrm{M}_{\Gamma}}
=\frac{\|\partial h^{(Q)}_k(\matr{I}_2)\|}{8\sqrt{5}\mathrm{M}_{\Gamma}}
\geq \frac{\varepsilon}{8\sqrt{5}\mathrm{M}_{\Gamma}}\|\matr{\Lambda}(\matr{X}_{k-1})\|. 
\end{align*}
Note that 
\begin{align*}
\|\matr{\Utwo}^{*}_{k}-\matr{I}_{2}\| = 2\sqrt{2}\left|\sin(\frac{\theta_{*}}{2})\right| \geq \frac{\sqrt{2}}{4}|\sin(4\theta_{*})|.
\end{align*}
We only need to set $\kappa_{Q}=\frac{\sqrt{2}\varepsilon}{32\sqrt{5}\mathrm{M}_{\Gamma}}$ in this case.
The proof is complete. 
\end{proof}

By \Cref{lemma:h_k_condition_SU-0}, \Cref{lemma:weak_rho_inequality_c-3}, inequalities \cref{eq:x_k_x_k-1_1} and \cref{eq:x_k_x_k-1_2}, we can now easily get the following result by setting $\eta_{Q}=(\kappa_{Q}\iota_{Q})/\mathrm{M}_{\omega}^2$. 

\begin{corollary}\label{lemma:weak_rho_inequality_c-4}
In \Cref{alg:BCD-G-RL} for cost function \cref{eq:sym_tensor_diagonalization-gener-0-equi},
if the iterates remain bounded, i.e., the condition \cref{eq:condi_bounded_X_gene} is satisfied, then 
\begin{equation*}
g(\matr{X}_{k-1})-g(\matr{X}_{k})\geq\eta_{Q}\|\ProjGrad{g}{\matr{X}_{k-1}}\|^2,
\end{equation*}
whenever the elementary function $h_k=h^{(Q)}_{k}$. 
\end{corollary}

\section{Convergence analysis}\label{sec:BCD-G}

In this section, based on the inequalities derived in \Cref{sec:trian_diago_rotation,sec:trian_diago_rotation-2}, we will prove our main results about the global and weak convergence of the BCD-G and Jacobi-G algorithms formulated in \Cref{subsecc:bcd-g}. 

\subsection{Convergence analysis of BCD-G algorithms}\label{subsec:convergence_BCD-GLU}

We now prove the following results about the global and weak convergence of BCD-G algorithms. 

\begin{theorem}\label{theorem_main_global_conv}
In BCD-GLU, BCD-GQU and BCD-GU algorithms for cost function \cref{eq:sym_tensor_diagonalization-gener-0-equi},
if the iterates remain bounded, i.e., the condition \cref{eq:condi_bounded_X_gene} is satisfied, 
then the iterates $\{\matr{\omega}_k\}_{k\geq 1}$ converge to a point $\matr{\omega}_{*}$.
\end{theorem}

\begin{proof}
We first prove the case of BCD-GLU algorithm. 
By \Cref{lemma:global_rho_c,lemma:global_rho_c-2}, we see that  
\begin{equation}\label{equation-condition-weak_c-5-0-proof-7.1}
g(\matr{X}_{k-1})-g(\matr{X}_{k})\geq\min(\sigma_{D},\sigma_{LU})\|\ProjGrad{g}{\matr{X}_{k-1}}\|\|\matr{X}_{k}-\matr{X}_{k-1}\|,
\end{equation}
whenever the elementary function $h_k=h^{(D)}_{k}, h^{(U)}_{k}$ or $h^{(L)}_k$. 
By the above inequality \cref{equation-condition-weak_c-5-0-proof-7.1} and \Cref{lemma:line_search}, we have that 
\begin{equation}\label{eq:sufficient_descent-2-proof-6.11}
f(\matr{\omega}_{k-1})-f(\matr{\omega}_{k})\geq\min(\sigma_{D},\sigma_{LU},\sigma_{p})\|\ProjGrad{f_{t_k}}{\matr{\omega}_{k-1}}\|\|\matr{\omega}_{k}-\matr{\omega}_{k-1}\|
\end{equation}
always holds in BCD-GLU algorithm, which is the inequality \cref{eq:sufficient_descent-2} in \Cref{theorem-SU15-2}, if we set $\sigma=\min(\sigma_{D},\sigma_{LU},\sigma_{p})$. 
Therefore, if $\matr{\omega}_*$ is an accumulation point of the iterates $\{\matr{\omega}_k\}_{k\geq 1}$ produced by BCD-GLU algorithms, it is the limit point.
Note that the iterates  $\{\matr{\omega}_k\}_{k\geq 1}$ remain bounded by condition \cref{eq:condi_bounded_X_gene}. There exists an accumulation point $\matr{\omega}_*$ such that the iterates $\{\matr{\omega}_k\}_{k\geq 1}$ converge to $\matr{\omega}_{*}$. 
For other two cases of BCD-GQU and BCD-GU algorithms, by \Cref{lemma:line_search}, \Cref{lemma:global_rho_c}, \Cref{lemma:global_rho_c-2} and \Cref{lemma:h_k_condition_SU}, we can similarly prove that the inequality \cref{eq:sufficient_descent-2} is always satisfied, if we set $\sigma=\min(\sigma_{D},\sigma_{LU},\sigma_{Q},\sigma_{p})$ and $\min(\sigma_{D},\sigma_{LU},\sigma_{p})$, respectively.  
The proof is complete. 
\end{proof}

To help the readers better understand the proof of \cref{theorem_main_global_conv}, we now summarize in \Cref{figure-flow-proof-3} the proof structure of \cref{theorem_main_global_conv} for BCD-GLU algorithm. Other two cases of BCD-GQU and BCD-GU algorithms are similar. 

\begin{figure}[h]
\centering
\begin{tikzpicture}[node distance=3.2cm]
\node (start1) [startstop] {\Cref{lemma:global_rho_c-2}};
\node (start2) [startstop, above of=start1, yshift=-2.2cm] {\Cref{lemma:global_rho_c-0}};
\node (start3) [startstop, right of=start1] {Inequality \cref{equation-condition-weak_c-5-0-proof-7.1}};
\node (start4) [startstop, right of=start3] {Inequality \cref{eq:sufficient_descent-2-proof-6.11}};
\node (start5) [startstop, right of=start4] {\cref{theorem_main_global_conv}};
\node (start6) [startstop, below of=start3, yshift=2.2cm]{\Cref{lemma:line_search}};
\node (start7) [startstop, above of=start4, yshift=-2.2cm] {\Cref{theorem-SU15-2}};
\node (start8) [startstop, above of=start3, yshift=-1.2cm] {\Cref{theorem-SU15}};
\node (start9) [startstop, above of=start5, yshift=-1.2cm] {Inequality \cref{eq:inequality_BCD_G}};
\node (start10) [startstop0, above of=start2, yshift=-2.2cm] {};
\node (start11) [startstop, below of=start1, yshift=2.2cm] {\Cref{lemma:global_rho_c-2-0}};
\node (start12) [startstop, above of=start3, yshift=-2.2cm] {\Cref{lemma:global_rho_c}};
\draw [arrow] (start1) -- (start3);
\draw [arrow] (start2) -- (start12);
\draw [arrow] (start3) -- (start4);
\draw [arrow] (start6) -| (start4);
\draw [arrow] (start4) -- (start5);
\draw [arrow] (start7) -| (start5);
\draw [arrow] (start8) -| (start7);
\draw [arrow] (start9) -| (start7);
\draw [arrow] (start11) -- (start1);
\draw [arrow] (start12) -- (start3);
\end{tikzpicture}
\caption{Proof structure of \cref{theorem_main_global_conv} for BCD-GLU algorithm.}
\label{figure-flow-proof-3}
\end{figure}

\begin{theorem}\label{theorem_main_weak_conv}
In BCD-GLU, BCD-GQU and BCD-GU algorithms for cost function \cref{eq:sym_tensor_diagonalization-gener-0-equi},
if the iterates remain bounded, i.e., the condition \cref{eq:condi_bounded_X_gene} is satisfied, 
and $\matr{\omega}_*$ is an accumulation point of the iterates $\{\matr{\omega}_{k}\}_{k\ge1}$,
then $\matr{\omega}_*$ is a stationary point of the cost function \cref{eq:sym_tensor_diagonalization-gener-0-equi}. 
\end{theorem}

\begin{proof}
We first prove the case of BCD-GLU algorithm. 
By \Cref{lemma:weak_rho_inequality_c}, we see that  
\begin{equation}\label{equation-condition-weak_c-5-0-proof-7.1-43}
g(\matr{X}_{k-1})-g(\matr{X}_{k})\geq\min(\eta_{D},\eta_{LU})\|\ProjGrad{g}{\matr{X}_{k-1}}\|^2,
\end{equation}
whenever the elementary function $h_k=h^{(D)}_{k}, h^{(U)}_{k}$ or $h^{(L)}_k$. 
By the above inequality \cref{equation-condition-weak_c-5-0-proof-7.1-43} and \Cref{lemma:line_search-2}, we have that 
\begin{equation}\label{eq:sufficient_descent-2-proof-6.11-32}
f(\matr{\omega}_{k-1})-f(\matr{\omega}_{k})\geq\min(\eta_{D},\eta_{LU},\eta_{p})\|\ProjGrad{f_{t_k}}{\matr{\omega}_{k-1}}\|^2
\end{equation}
always holds in BCD-GLU algorithm, which is the inequality \cref{equation-condition-weak_c} in \Cref{theorem-weak-general_c}, if we set $\eta=\min(\eta_{D},\eta_{LU},\eta_{p})$.
Therefore, if $\matr{\omega}_*$ is an accumulation point of the iterates $\{\matr{\omega}_k\}_{k\geq 1}$ produced by BCD-GLU algorithm,
then $\matr{\omega}_*$ is a stationary point. 
For other two cases of BCD-GQU and BCD-GU algorithms, by \Cref{lemma:weak_rho_inequality_c} and \Cref{lemma:weak_rho_inequality_c-4},
we can similarly prove that the inequality \cref{equation-condition-weak_c} is always satisfied, if we set $\eta=\min(\eta_{D},\eta_{LU},\eta_{Q},\eta_{p})$ and $\min(\eta_{D},\eta_{LU},\eta_{p})$, respectively.  
The proof is complete. 
\end{proof}

To help the readers better understand the proof of \cref{theorem_main_weak_conv}, we now summarize in \Cref{figure-flow-proof-4-5} the proof structure of \cref{theorem_main_weak_conv} for BCD-GLU algorithm. Other two cases of BCD-GQU and BCD-GU algorithms are similar. 

\begin{figure}[h]
\centering
\begin{tikzpicture}[node distance=3.2cm]
\node (start1) [startstop] {\Cref{lemma:global_rho_c-2-0}};
\node (start2) [startstop0, above of=start1, yshift=-2.2cm] {};
\node (start3) [startstop, right of=start1] {Inequality \cref{equation-condition-weak_c-5-0-proof-7.1-43}};
\node (start4) [startstop, right of=start3] {Inequality \cref{eq:sufficient_descent-2-proof-6.11-32}};
\node (start5) [startstop, right of=start4] {\cref{theorem_main_weak_conv}};
\node (start7) [startstop, above of=start5, yshift=-2.2cm] {\Cref{theorem-weak-general_c}};
\node (start12) [startstop, above of=start3, yshift=-2.2cm]{\Cref{lemma:weak_rho_inequality_c}};
\node (start8) [startstop, above of=start12, yshift=-2.2cm] {\Cref{lemma:weak_rho_inequality_c-1}};
\node (start10) [startstop, above of=start2, yshift=-2.2cm] {\Cref{lemma:global_rho_c-0}};
\node (start13) [startstop, above of=start4, yshift=-2.2cm] {\Cref{lemma:line_search-2}};
\node (start9) [startstop, above of=start13, yshift=-2.2cm] {Inequality \cref{eq:inequality_BCD_G}};
\draw [arrow] (start1) |- (start12);
\draw [arrow] (start3) -- (start4);
\draw [arrow] (start4) -- (start5);
\draw [arrow] (start7) -- (start5);
\draw [arrow] (start8) -- (start12);
\draw [arrow] (start9) -| (start7);
\draw [arrow] (start10) |- (start12);
\draw [arrow] (start12) -- (start3);
\draw [arrow] (start13) -- (start4);
\end{tikzpicture}
\caption{Proof structure of \cref{theorem_main_weak_conv} for BCD-GLU algorithm.}
\label{figure-flow-proof-4-5}
\end{figure}

\subsection{Convergence analysis of Jacobi-G algorithms}
 
Similar as in \Cref{subsec:convergence_BCD-GLU}, we have the following results about the global and weak convergence of Jacobi-G algorithms. We omit the detailed proofs here. 

\begin{theorem}\label{theorem:weak_Jaco_Lu-0}
In Jacobi-GLU and Jacobi-GQU algorithms for cost function \cref{eq:sym_tensor_diagonalization-gener-1},
if the iterates remain bounded, i.e., the condition \cref{eq:condi_bounded_X_gene-2} is satisfied, 
then the iterates $\{\matr{X}_{k}\}_{k\ge1}$ converge to a point $\matr{X}_{*}$. 
\end{theorem}

\begin{theorem}\label{theorem_main_weak_conv-Jacobi-GLU-09}
In Jacobi-GLU and Jacobi-GQU algorithms for cost function \cref{eq:sym_tensor_diagonalization-gener-1},
if the iterates remain bounded, i.e., the condition \cref{eq:condi_bounded_X_gene-2} is satisfied, 
and $\matr{X}_*$ is an accumulation point of the iterates $\{\matr{X}_{k}\}_{k\ge1}$,
then $\matr{X}_*$ is a stationary point of the cost function \cref{eq:sym_tensor_diagonalization-gener-1}.
\end{theorem}

\begin{remark}\label{subsec:Glu-M}
We propose two natural variants of Jacobi-GLU and Jacobi-GQU algorithms, which will be called \emph{Jacobi-GLU-M} and \emph{Jacobi-GQU-M} algorithms, respectively.
In these two algorithms, in each iteration, among all the index pairs $(i_k,j_k)$ and elementary functions $h_k$ satisfying inequality \cref{eq-inequality-G-LU_c}, 
we choose $(i_k,j_k)$ and $h_k$ such that the cost function obtains the largest reduction.
It is clear that \Cref{theorem:weak_Jaco_Lu-0,theorem_main_weak_conv-Jacobi-GLU-09} also apply to these two new variants. 
\end{remark}

\begin{remark}\label{remar:jacobi-c}
In Jacobi-GLU and Jacobi-GQU algorithms, a more natural way of choosing the index pair $(i_k,j_k)$ is according to a \emph{cyclic} ordering.
In fact, this cyclic way has often been used in the literature \cite{maurandi2014decoupled,sorensen2009approximate,wang2012complex}. 
In this case,
we call them the \emph{Jacobi-CLU} and \emph{Jacobi-CQU} algorithms, respectively.  
\end{remark}

\section{Numerical experiments}\label{sec-experiment-hooi}
In the BCD-G and Jacobi-G algorithms of this paper, there exist several parameters to be adjusted, including the positive constant $\upsilon$ in inequality \cref{eq:inequality_BCD_G},  the stepsize $t_{k-1}$ in equation \cref{eq:U_k_update}, the positive constant $\varepsilon$ in inequality \cref{eq-inequality-G-LU_c},  $\varsigma_{D}>0$ in \cref{remar:rho_setting_four_case_c}, and $\varsigma_{Q}>0$ in  \cref{remark:condition_inner}. 
In this section, we choose different values for the positive constant $\varepsilon$ in inequality \cref{eq-inequality-G-LU_c}, while fixing other parameters as small positive constants. 
We set $(\cdot)^\TH = (\cdot)^\H$ in both the cost functions \cref{eq:sym_tensor_diagonalization-gener-0-equi} and \cref{eq:sym_tensor_diagonalization-gener-1}. 
All the algorithms run at most 1000 iterations. All the randomly generated complex matrices are uniformly distributed. 
All the computations are done using MATLAB R2019a. 
The numerical experiments are conducted on a PC with an $\text{Intel}^{\textregistered}$ $\text{Core}^{\text{TM}}$ i5 CPU at 2.11 GHz and 8.00 GB of RAM in 64bt Windows operation system.

\begin{figure}%[htb!]
\centering
\begin{minipage}[b]{0.44\linewidth}
%\centering
{\includegraphics[height=5.2cm]{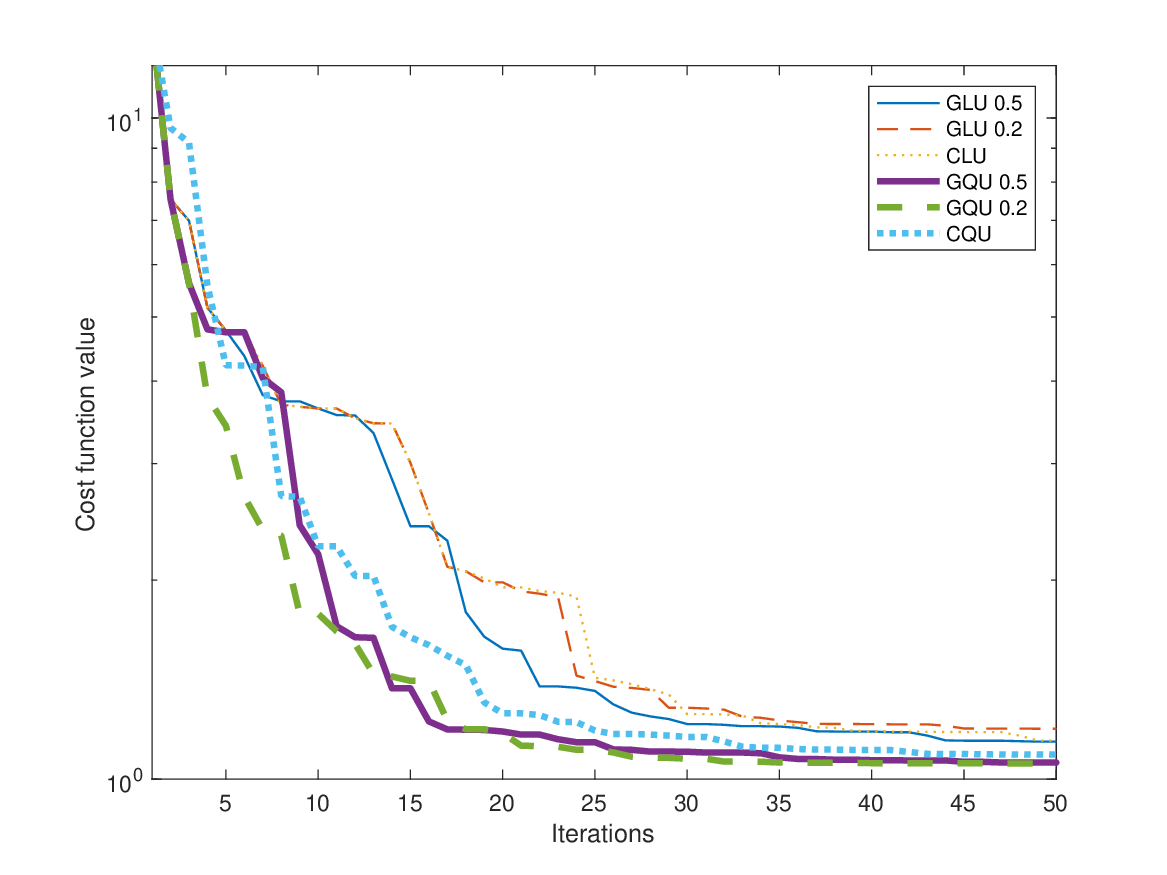}}
% \vspace{1.5cm}
\centerline{\quad \quad \quad \quad\quad \quad \quad (i)}\medskip
\end{minipage}
\hfill
\begin{minipage}[b]{0.44\linewidth}
%	\centering
{\includegraphics[height=5.2cm]{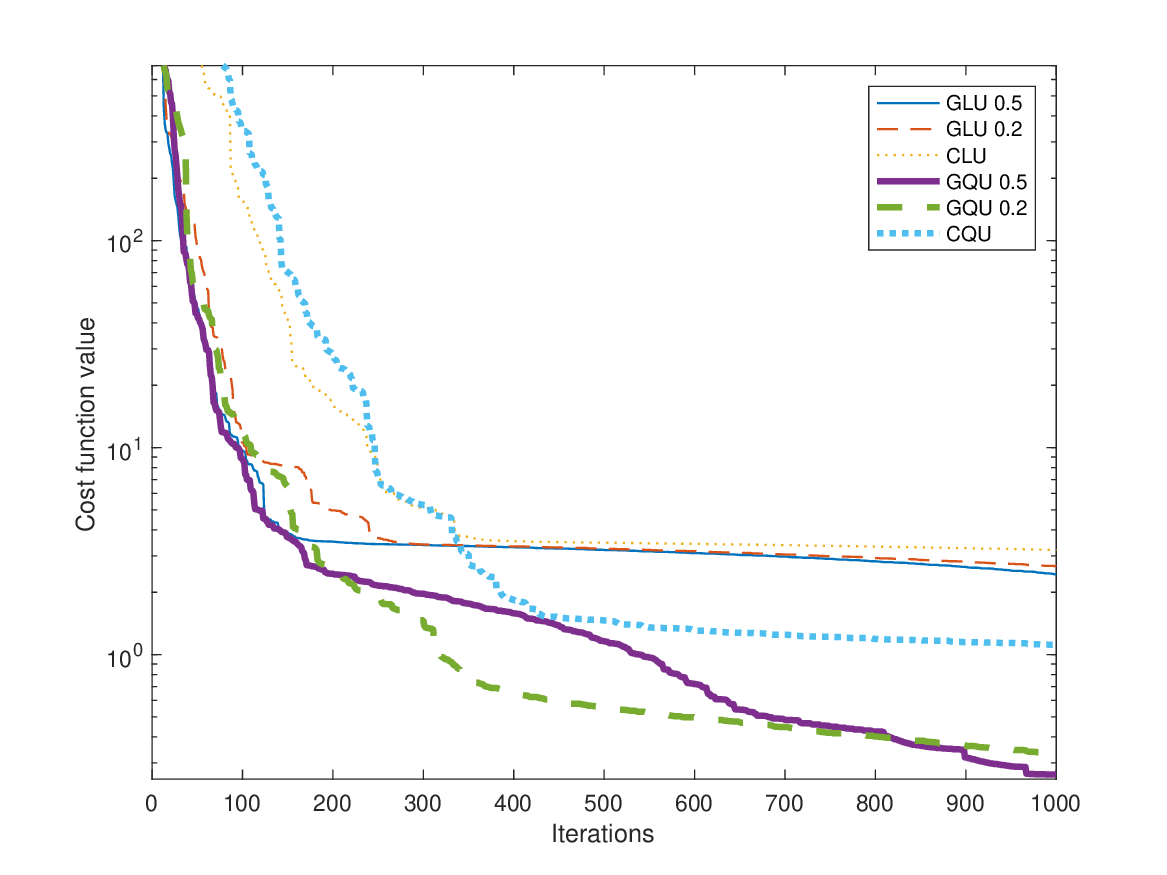}}
% \vspace{1.5cm}
\centerline{\quad \quad \quad \quad\quad \quad \quad (ii)}\medskip\end{minipage}
\hfill
\begin{minipage}[b]{0.44\linewidth}
%	\centering
{\includegraphics[height=5.2cm]{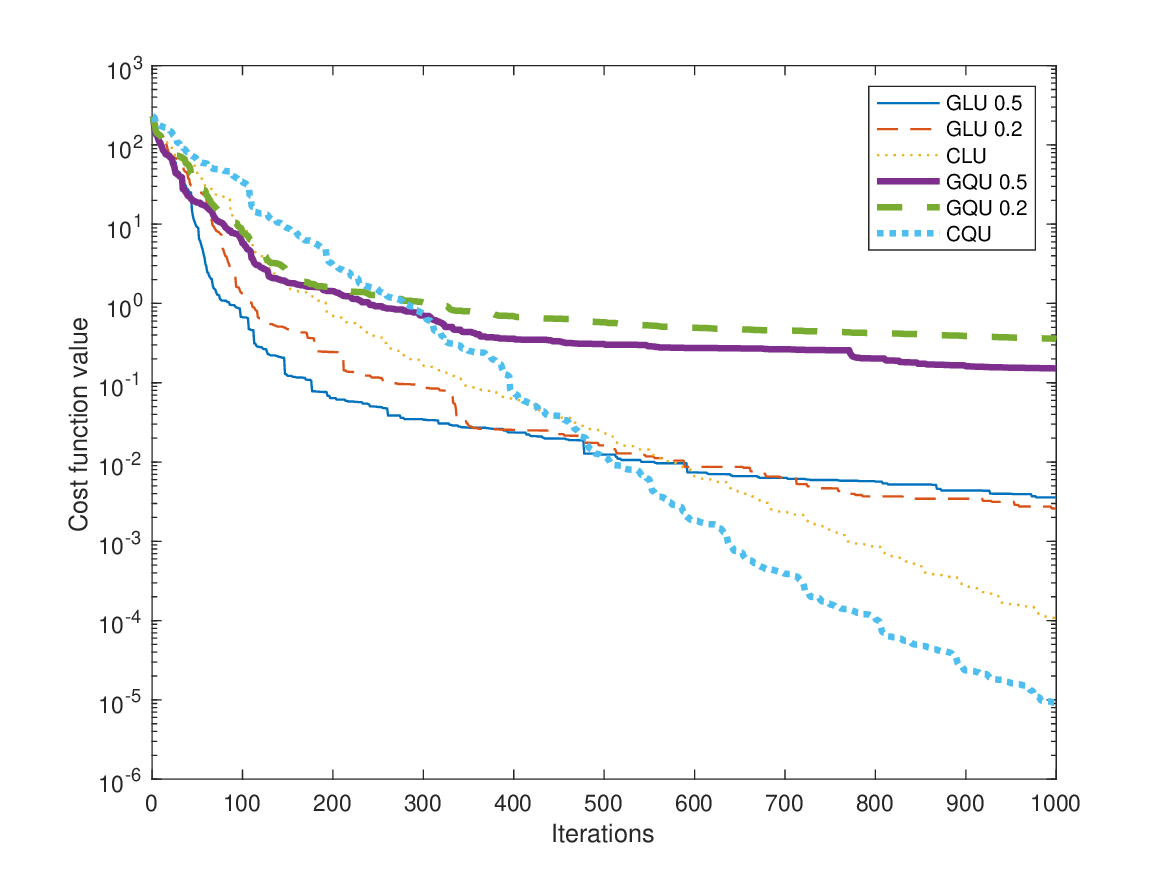}}
% \vspace{1.5cm}
\centerline{\quad \quad \quad \quad\quad \quad \quad (iii)}\medskip\end{minipage}
\hfill
\begin{minipage}[b]{0.44\linewidth}
%	\centering
{\includegraphics[height=5.2cm]{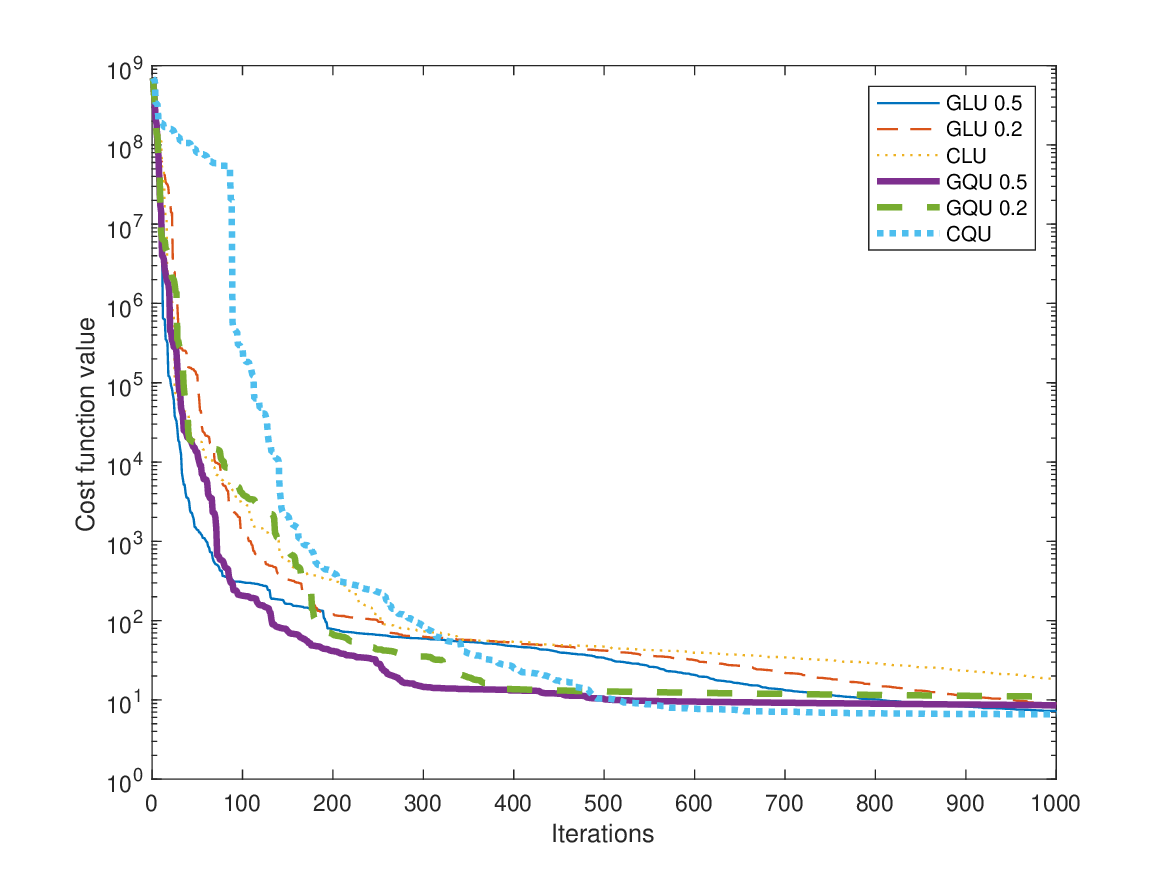}}
% \vspace{1.5cm}
\centerline{\quad \quad \quad \quad\quad \quad \quad (iv)}\medskip\end{minipage}
\caption{Experimental results for BCD-G algorithms in \Cref{example-1-1}.}
\label{figure-example-1-1}
\end{figure}

\begin{example1}\label{example-1-1}
For the following sets of complex matrices, 
we run BCD-GLU and BCD-GQU algorithms to minimize the cost function \cref{eq:sym_tensor_diagonalization-gener-0-equi}. 
The values of cost function \cref{eq:sym_tensor_diagonalization-gener-0-equi} in the iterations are shown in \Cref{figure-example-1-1}. The positive constant $\upsilon$ in inequality \cref{eq:inequality_BCD_G} is fixed to $0.001$. 
For the positive constant $\varepsilon$ in inequality \cref{eq-inequality-G-LU_c}, we choose different values.  
For example, \emph{BCD-GLU 0.5} means the BCD-GLU algorithm with  $\varepsilon=0.5\sqrt{\frac{2}{3m(m-1)}}$. 
If $\varepsilon=0$, we denote the BCD-GLU and BCD-GQU algorithms by BCD-CLU and BCD-CQU, respectively. The starting point is $\matr{\omega}_{0}=(\matr{I}_{n\times m},\matr{I}_{m})$. \\
(i) We set $n=5$, $m=3$, and randomly generate complex matrices $\{\matr{A}_{\ell}\}_{1\leq\ell\leq 3}\subseteq\CC^{5\times 5}$.\\
(ii) We set $n=10$, $m=8$, randomly generate a complex matrix $\matr{X}\in\CC^{10\times 10}$, and set 
$\matr{A}^{(\ell)}=\matr{X}^\H (\matr{I}_{10}+\vect{e}_{\ell}^{\T}\vect{e}_{\ell})\matr{X}$
for $1\leq \ell\leq 5$.\\
(iii) We set $n=10$, $m=8$, randomly generate a complex upper triangular matrix $\matr{X}\in\UT_{10}(\CC)$, complex diagonal matrices  $\{\matr{D}_{\ell}\}_{1\leq\ell\leq 5}\subseteq\CC^{10\times 10}$, and set $\matr{A}^{(\ell)}=\matr{X}^\H \matr{D}_{\ell} \matr{X}$ for $1\leq \ell\leq 5$. \\
(iv) We set $n=10$, $m=8$, randomly generate a complex nonsingular matrix $\matr{X}\in\SL_{10}(\CC)$, complex diagonal matrices  $\{\matr{D}_{\ell}\}_{1\leq\ell\leq 5}\subseteq\CC^{10\times 10}$, and set $\matr{A}^{(\ell)}=\matr{X}^\H \matr{D}_{\ell} \matr{X}$ for $1\leq \ell\leq 5$.
\end{example1} 

\begin{figure}%[htb!]
\centering
\begin{minipage}[b]{0.44\linewidth}
%\centering
{\includegraphics[height=5.2cm]{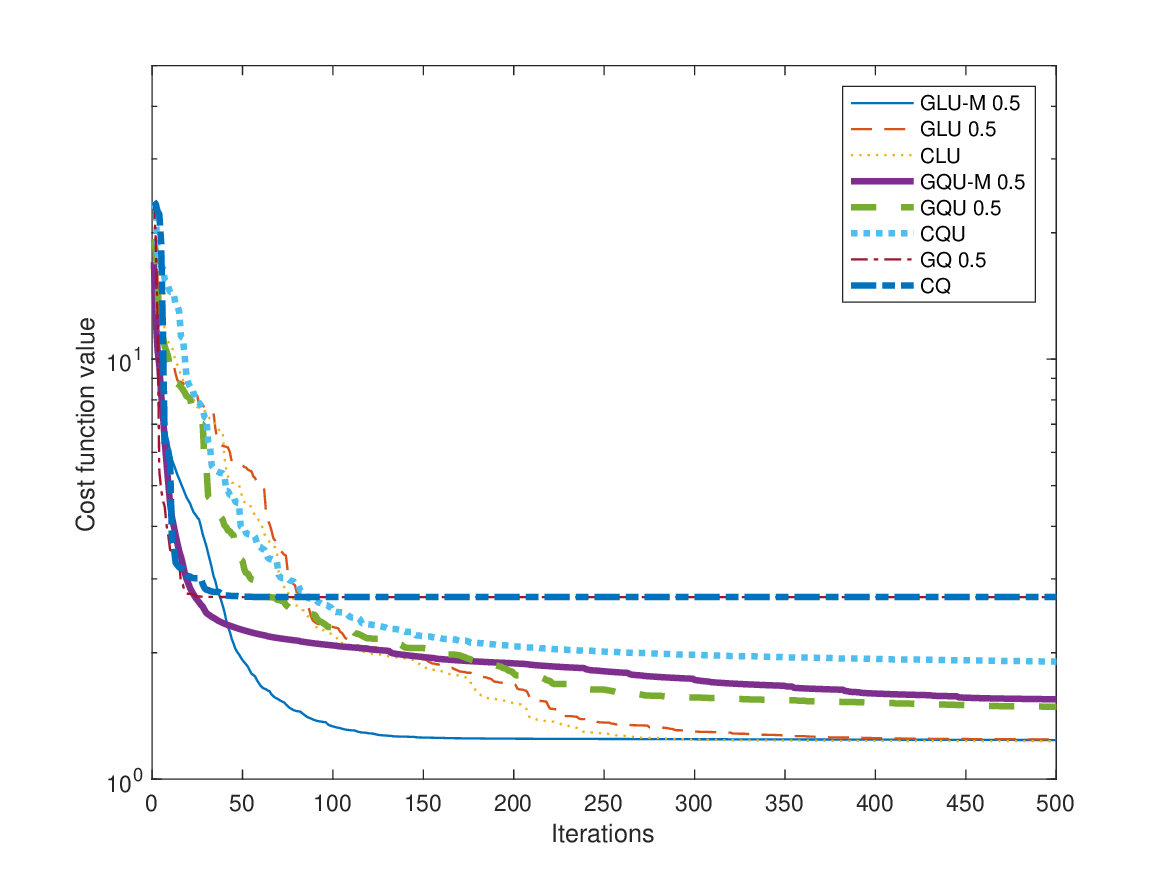}}
% \vspace{1.5cm}
\centerline{\quad \quad \quad \quad\quad \quad \quad (i)}\medskip
\end{minipage}
\hfill
\begin{minipage}[b]{0.44\linewidth}
%	\centering
{\includegraphics[height=5.2cm]{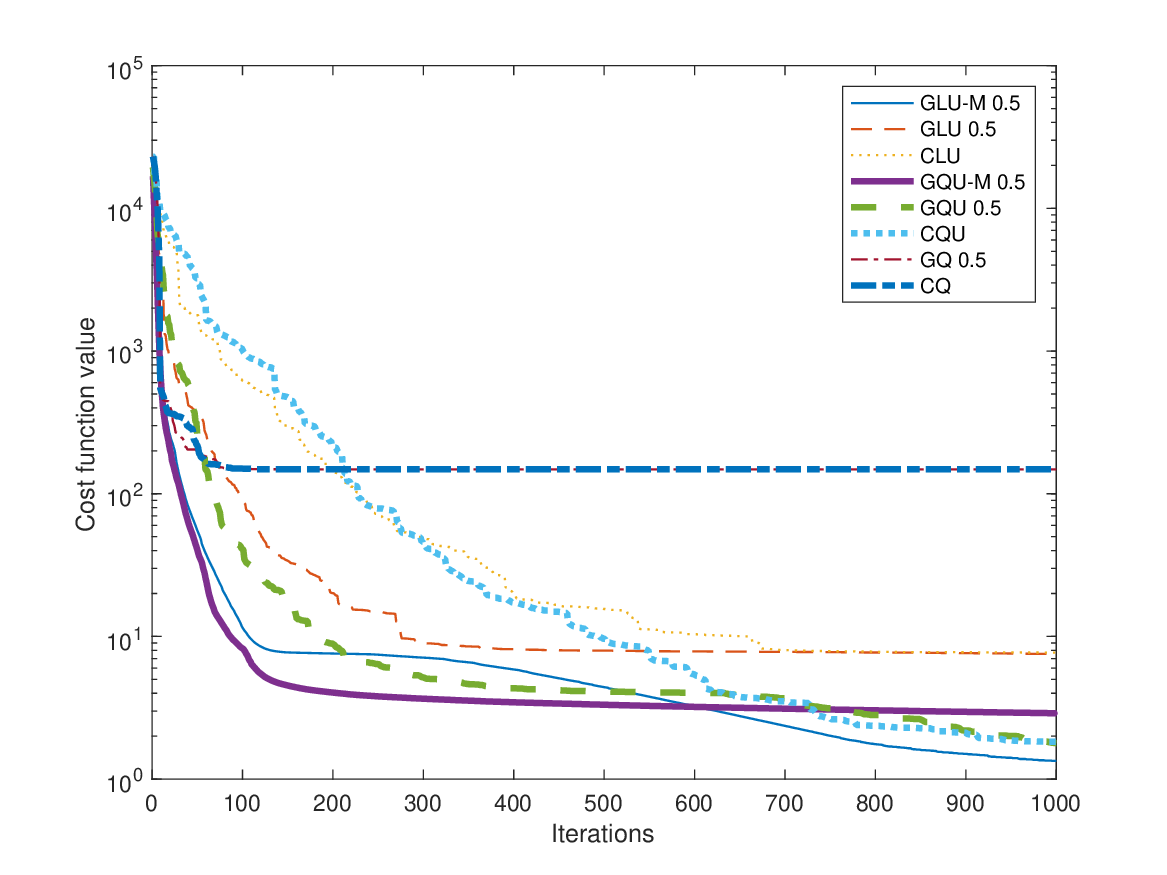}}
% \vspace{1.5cm}
\centerline{\quad \quad \quad \quad\quad \quad \quad (ii)}\medskip\end{minipage}
\caption{Experimental results for Jacobi-type algorithms in \Cref{example-1}.}
\label{figure-example-1}
\end{figure}

\begin{example1}\label{example-1}
For the following sets of complex matrices, 
we run eight Jacobi-type algorithms to minimize the cost function \cref{eq:sym_tensor_diagonalization-gener-1}.
Here, we denote by Jacobi-GQ the gradient-based Jacobi-type algorithm on the unitary group proposed in \cite{ULC2019}, and by Jacobi-CQ the Jacobi-type algorithm on the unitary group with a cyclic ordering. 
Note that Jacobi-GQ and Jacobi-CQ find the iterates only in $\UNN{m}(\CC)$, not in $\SL_m(\CC)$.
The values of cost function \cref{eq:sym_tensor_diagonalization-gener-1} in the iterations are shown in \Cref{figure-example-1}.  
We choose the starting point $\matr{X}_{0}=\matr{I}_{m}$.\\ 
(i) We randomly generate two complex matrices $\{\matr{A}_{\ell}\}_{1\leq\ell\leq 2}\subseteq\CC^{5\times 5}$.\\
(ii) We randomly generate a complex matrix $\matr{X}\in\CC^{10\times 10}$, and set 
$\matr{A}^{(\ell)}=\matr{X}^\H (\matr{I}_{10}+\vect{e}_{\ell}^{\T}\vect{e}_{\ell})\matr{X}$
for $1\leq \ell\leq 10$. 
\end{example1}

From the above numerical experiments, we can see that: (i) in \Cref{example-1-1}, compared with BCD-GLU algorithms, the BCD-GQU algorithms generally have better performances;  
(ii) in \Cref{example-1}, compared with Jacobi-GQU algorithms, the Jacobi-GLU algorithms generally have better performances; 
(iii) in \Cref{example-1}, compared with the Jacobi-GQ and Jacobi-CQ algorithms on $\UNN{m}(\CC)$, the Jacobi-type algorithms on $\SL_{m}(\CC)$ considered in this paper always obtain much smaller cost function values, and they also need more iterations to attain steady state values of cost functions. 

\section{Conclusions}\label{sec:conclusion}

In this paper, to solve JADM problem \cref{eq:sym_tensor_diagonalization-gener-0}, which is important in BSS problem, we formulate two different equivalent formulations, \emph{i.e.}, problem \cref{eq:sym_tensor_diagonalization-gener-0-equi} defined on $\St(m,n,\CC)\times\SL_{m}(\CC)$, and problem \cref{eq:sym_tensor_diagonalization-gener-1} defined on $\SL_{m}(\CC)$. Then, for these two approaches, based on the Riemannian gradients, we propose three BCD-G algorithms and two Jacobi-G algorithms, and establish their global and weak convergence, under the condition that the iterates are bounded. 
An interesting question is, in the BCD-G and Jacobi-G algorithms, whether one can find a method to guarantee both the boundedness of the iterates and the inequalities \cref{eq:sufficient_descent-2,equation-condition-weak_c} for global and weak convergence. 
If so, then one can get rid of the dependence of convergence results on the condition that the iterates are bounded.  

\appendix

\section{Proofs in \Cref{sec:Jacobi_g_rota}}\label{sec:ml_proofs_wellness_Jaco_part_01} Before the proof of \cref{lem:ProjGradSubmatrix}, we need to show a lemma, which is similar as equation \cref{eq:gradient_exp_map} and can be directly obtained from \cite[Eq. (3.31)]{absil2009optimization}. 

\begin{lemma}\label{lem:gra_rela_elem}
Let $\exp: \CC^{2\times 2}  \to \GL_{2}(\CC)$ 
be the matrix exponential function \cite{absil2009optimization,baker2012matrix,Hall15:lie} sending $\matr{\Delta}$ to
$\exp(\matr{\Delta})$. \\
(i) If $h:\SUNN{2}(\CC)\to\RR$ is a differentiable function and $\matr{\Delta}\in\su_2(\CC)= \Tang{\matr{I}_2}{\SUNN{2}(\CC)}$, we have that 
\begin{equation}\label{eq:gradient_exp_map-2}
\langle\matr{\Delta},\ProjGrad{h}{\matr{I}_2}\rangle_{\matr{I}_2}  =
\left.\left( \frac{d}{dt} h (\exp(t\matr{\Delta}) )\right)\right|_{t=0}. 
\end{equation}
(ii) If $h:\SUT_2(\CC)\to\RR$ is a differentiable function and $\matr{\Delta}\in\sut_2(\CC)= \Tang{\matr{I}_2}{\SUT_2(\CC)}$, we have the relationship \cref{eq:gradient_exp_map-2}. \\ 
(iii) If $h:\SLT_2(\CC)\to\RR$ is a differentiable function and $\matr{\Delta}\in\slt_2(\CC)= \Tang{\matr{I}_2}{\SLT_2(\CC)}$, we have the relationship \cref{eq:gradient_exp_map-2}.\\ 
(iv) If $h:\DT_2(\CC)\to\RR$ is a differentiable function and $\matr{\Delta}\in\sdt_2(\CC)= \Tang{\matr{I}_2}{\DT_2(\CC)}$, we have the relationship \cref{eq:gradient_exp_map-2}. 
\end{lemma}

\begin{proof}[Proof of \cref{lem:ProjGradSubmatrix}]
Define a projection operator $\mathcal{P}_{i,j}: \CC^{m\times m} \to\CC^{2\times 2}$ extracting a submatrix of $\matr{X} \in \CC^{m\times m}$ as in \cite[Eq. (3.7)]{ULC2019}, and $\mathcal{P}_{i,j}^{\T}: \CC^{2\times 2} \to\CC^{m\times m}$ the conjugate operator.  
For the elementary function $\hij{i}{j}{\matr{X}}$ defined in \cref{eq:elemn_func_h}, if $\matr{\Delta} \in \su_2(\CC) = \Tang{\matr{I}_2}{\SUNN{2}(\CC)}$, we have that 
\begin{align*}
&\langle\matr{\Delta},\ProjGrad{\hij{i}{j}{\matr{X}}}{\matr{I}_2}\rangle_{\matr{I}_2}  =
\left.\left( \frac{d}{dt} \hij{i}{j}{\matr{X}} (\exp(t\matr{\Delta}) )\right)\right|_{t=0} \ \ (\text{by \Cref{lem:gra_rela_elem}(i)}) \\
&= 
\left.\left( \frac{d}{dt} g (\matr{X} \Gmat{i}{j}{\exp(t\matr{\Delta} )} )\right)\right|_{t=0}
= \left.\left( \frac{d}{dt} g (\text{Exp}_{\matr{X}} (\matr{X}\mathcal{P}_{i,j}^{\T}(\matr{\Delta} )t )) \right)\right|_{t=0} \\
&= \langle\matr{X} \mathcal{P}_{i,j}^{\T}(\matr{\Delta}),\ProjGrad{g}{\matr{X}}\rangle_{\matr{X}} \ \ (\text{by equation}\   \cref{eq:gradient_exp_map})=
\langle\matr{\Delta},\mathcal{P}_{i,j} (\matr{\Lambda})\rangle_{\matr{I}_2},
\end{align*}
where $\text{Exp}_{\matr{X}}$ is the map defined in \cref{eq:Exp_x_sl}. 
Note that $\matr{\Delta} \in \su_2(\CC)$ and $\ProjGrad{\hij{i}{j}{\matr{X}}}{\matr{I}_2} \in \su_2(\CC)$. 
The result can be obtained by direct calculations. 
For other three elementary functions  $h^{(U)}_{(i,j),\matr{X}}$, $h^{(L)}_{(i,j),\matr{X}}$ and $h^{(D)}_{(i,j),\matr{X}}$, similar as the above case, we can obtain the results by  \Cref{lem:gra_rela_elem}(ii), \Cref{lem:gra_rela_elem}(iii) and \Cref{lem:gra_rela_elem}(iv), respectively. 
The proof is complete.
\end{proof}

We need a simple lemma before the proofs of \Cref{theorem-wellness-Jacobi-G-LU_c} and \Cref{theorem-wellness-Jacobi-G-GU_c}.

\begin{lemma}\label{lemma-ai-sum-0_c}
(i) If $z_1,z_2\in\CC$, 
then 
\begin{equation*}
|z_1-z_2|^2 + |z_2|^2 \geq \frac{3-\sqrt{5}}{2}(|z_1|^2+|z_2|^2).
\end{equation*}
(ii) If $\{z_i\}_{1\leq i\leq m}\subseteq\CC$ satisfy $\sum_{1\leq i\leq m}z_i=0$, 
then 
\begin{equation*}
\sum_{1\leq i<j\leq m}|z_i-z_j|^2 = m\sum_{1\leq i\leq m}|z_i|^2.
\end{equation*}
\end{lemma}

\begin{proof}[Proof of \Cref{theorem-wellness-Jacobi-G-LU_c}]
(i) We first prove the existence of such an index pair $(i_k,j_k)$ and an elementary function $h_{k}$ in \autoref{alg:jacobi-LU-G_c}. 
By \Cref{lem:ProjGradSubmatrix_2}, \Cref{lemma-ai-sum-0_c} and $\matr{\Lambda}=\matr{\Lambda}(\matr{X}_{k-1})\in\sll_{m}(\CC)$, we have that 
\begin{align*}
&\sum_{1\leq i_k<j_k\leq m}  \left(\|\partial h^{(U)}_k(\matr{I}_2)\|^2+\|\partial h^{(L)}_k(\matr{I}_2)\|^2+\|\partial h^{(D)}_k(\matr{I}_2)\|^2\right) \\
&= \sum_{1\leq i_k<j_k\leq m} \left(|\Lambda_{i_kj_k}|^2 + |\Lambda_{j_ki_k}|^2\right) + m\sum_{1\leq i_k\leq m} |\Lambda_{i_ki_k}|^2 \geq \|\matr{\Lambda}\|^2. 
\end{align*}
Therefore, there exist an index pair $(i_k,j_k)$ and an elementary function $h_k=h^{(U)}_k, h^{(L)}_k$ or $h^{(D)}_k$ such that $\frac{3}{2}m(m-1)\|\partial h_k(\matr{I}_2)\|^2\geq\|\matr{\Lambda}\|^2$. \\
(ii) We now prove the existence in \autoref{alg:jacobi-LQ-G_c}. 
Similar as above, 
we get that 
\begin{align*}
&\sum_{1\leq i_k<j_k\leq m}  \left(\|\partial h^{(Q)}_k(\matr{I}_2)\|^2+\|\partial h^{(U)}_k(\matr{I}_2)\|^2 + \|\partial h^{(D)}_k(\matr{I}_2) \|^2\right) \\
&= \sum_{1\leq i_k<j_k\leq m}  \left(|\Lambda_{i_kj_k}^{*} - \Lambda_{j_ki_k}|^2+|\Lambda_{i_kj_k}|^2 + |\Lambda_{i_ki_k} - \Lambda_{j_kj_k}|^2\right)\\
&\geq \frac{3-\sqrt{5}}{2}\sum_{1\leq i_k<j_k\leq m} \left(|\Lambda_{i_kj_k}|^2 + |\Lambda_{j_ki_k}|^2\right) + m\sum_{1\leq i_k\leq m} |\Lambda_{i_ki_k}|^2
\geq \frac{3-\sqrt{5}}{2}\|\matr{\Lambda}\|^2. 
\end{align*}
Therefore, there exists an index pair $(i_k,j_k)$ and an elementary function $h_k =h^{(Q)}_k, h^{(U)}_k$ or $h^{(D)}_k$ such that
$\frac{3}{3-\sqrt{5}}m(m-1)\|\partial h_k(\matr{I}_2)\|^2\geq\|\matr{\Lambda}\|^2.$
The proof is complete. 
\end{proof}

\begin{proof}[Proof of \cref{theorem-wellness-Jacobi-G-GU_c}]
Note that the starting point  $\matr{X}_{0}\in\EUT_{m}(\CC)$ in \autoref{alg:jacob-GU-G_c}. 
We see that $\matr{X}_{k}\in\EUT_{m}(\CC)$ for all $k\in\NN$. 
By \Cref{lem:ProjGradSubmatrix_2}, \Cref{lemma-ai-sum-0_c} and $\matr{\Lambda}=\matr{\Lambda}(\matr{X}_{k-1})\in\eut_{m}(\CC)$, we get that 
\begin{align*}{\small
\sum_{1\leq i_k<j_k\leq m}  \left(\|\partial h^{(U)}_k(\matr{I}_2)\|^2 + \|\partial h^{(D)}_k(\matr{I}_2) \|^2\right)}
&{\small= \sum_{1\leq i_k<j_k\leq m}  \left(|\Lambda_{i_kj_k}|^2 + |\Lambda_{i_ki_k} - \Lambda_{j_kj_k}|^2\right)}\\
&{\small = \sum_{1\leq i_k<j_k\leq m} |\Lambda_{i_kj_k}|^2 + m\sum_{1\leq i_k\leq m} |\Lambda_{i_ki_k}|^2
\geq \|\matr{\Lambda}\|^2.}
\end{align*}
Therefore, there exist an index pair $(i_k,j_k)$ and an elementary function $h_k = h^{(U)}_k$ or $h^{(D)}_k$ such that 
$m(m-1)\|\partial h_k(\matr{I}_2)\|^2\geq\|\matr{\Lambda}\|^2$. 
The proof is complete. 
\end{proof}
\section{Proofs in \Cref{sec:trian_diago_rotation}}\label{sec:proofs_section_5}

\begin{proof}[Proof of \cref{lemma:global_rho_c-0}]
We now prove the inequality \cref{equation-condition-weak_c-5} by \Cref{lemm:elemen_express_c}(iii) in three different cases shown in \Cref{remar:rho_setting_four_case_c}.
\begin{itemize}
\item 
If $\gamma_{1}=\gamma_{2} = 0$, it is clear that the inequality \cref{equation-condition-weak_c-5} is satisfied for any $\iota_{D}>0$.\\
\item 
If $\varpi\in[0,\varsigma_{D})$, we get that 
\begin{align*}
&h^{(D)}_{k}(1,0) - h^{(D)}_{k}(x_k^*,y_k^*)
=\frac{3}{4}\gamma_{1}(1-4\varpi)  
=\frac{3(1-4\varpi)}{8(1-\varpi)}|\partial h^{(D)}_{k}(1,0)|\\
&\geq \frac{3(1-4\varpi)\varepsilon}{8(1-\varpi)}\|\matr{\Lambda}(\matr{X}_{k-1})\|
=\frac{3(1-4\varpi)\varepsilon}{4\sqrt{5}(1-\varpi)}\frac{\sqrt{5}}{2}\|\matr{\Lambda}(\matr{X}_{k-1})\|\\
&\geq\frac{3(1-4\varpi)\varepsilon}{4\sqrt{5}(1-\varpi)}\|\matr{\Lambda}(\matr{X}_{k-1})\|\|\matr{\Utwo}^{*}_{k}-\matr{I}_{2}\| 
\geq\frac{3(1-4\varsigma_{D})\varepsilon}{4\sqrt{5}}\|\matr{\Lambda}(\matr{X}_{k-1})\|\|\matr{\Utwo}^{*}_{k}-\matr{I}_{2}\|. 
\end{align*}
\item If $\varpi\in(\frac{1}{\varsigma_{D}},+\infty]$, we similarly get the above inequality. 
% \begin{align*}
% &h^{(D)}_{k}(1,0) - h^{(D)}_{k}(x_k^*,y_k^*)
% =
% \frac{3}{4}\gamma_{2}(1-\frac{4}{\varpi}) 
% =\frac{3(1-\frac{4}{\varpi})}{8(1-\frac{1}{\varpi})}|\partial h^{(D)}_{k}(1,0)|\\&
% \geq \frac{3(1-\frac{4}{\varpi})\varepsilon}{8(1-\frac{1}{\varpi})}\|\matr{\Lambda}(\matr{X}_{k-1})\|
% =\frac{3(1-\frac{4}{\varpi})\varepsilon}{4\sqrt{5}(1-\frac{1}{\varpi})}\frac{\sqrt{5}}{2}\|\matr{\Lambda}(\matr{X}_{k-1})\|\\
% &\geq\frac{3(1-\frac{4}{\varpi})\varepsilon}{4\sqrt{5}(1-\frac{1}{\varpi})}\|\matr{\Lambda}(\matr{X}_{k-1})\|\|\matr{\Utwo}^{*}_{k}-\matr{I}_{2}\|
% \geq\frac{3(1-4\varsigma_{D})\varepsilon}{4\sqrt{5}}\|\matr{\Lambda}(\matr{X}_{k-1})\|\|\matr{\Utwo}^{*}_{k}-\matr{I}_{2}\|. 
% \end{align*}
\item If $\varpi\in [\varsigma_{D},\frac{1}{\varsigma_{D}}]$,
it is easy to verify that 
\begin{equation}\label{eq:inequa_rho_r4_c}
\sqrt[\leftroot{-3}\uproot{5}4]{\gamma_{1}\gamma_{2}}\geq \frac{\sqrt[\leftroot{-2}\uproot{4}4]{\varsigma_{D}}}{2} \left(\sqrt{\gamma_{1}}+\sqrt{\gamma_{2}}\right).
\end{equation}
Then, we get that 
\begin{align*}
&h^{(D)}_{k}(1,0) - h^{(D)}_{k}(x_k^*,y_k^*)
= \left(\sqrt{\gamma_{1}}-\sqrt{\gamma_{2}}\right)^2
=\frac{1}{2}|\partial h^{(D)}_{k}(1,0)|\frac{\left|\sqrt{\gamma_{1}}-\sqrt{\gamma_{2}}\right|}{\sqrt{\gamma_{1}}+\sqrt{\gamma_{2}}}\\
&\geq \frac{\varepsilon\sqrt[\leftroot{-2}\uproot{4}4]{\varsigma_{D}}}{4}\|\matr{\Lambda}(\matr{X}_{k-1})\|\frac{\left|\sqrt{\gamma_{1}}-\sqrt{\gamma_{2}}\right|}{\sqrt[\leftroot{-3}\uproot{5}4]{\gamma_{1}\gamma_{2}}} \ \ (\textrm{by equation \cref{eq:inequa_rho_r4_c}})\\
&\geq\frac{\varepsilon\sqrt[\leftroot{-2}\uproot{4}4]{\varsigma_{D}}}{4}\|\matr{\Lambda}(\matr{X}_{k-1})\|
\frac{\left|\sqrt[\leftroot{-3}\uproot{5}4]{\gamma_{1}}-\sqrt[\leftroot{-3}\uproot{5}4]{\gamma_{2}}\right|\left(\sqrt{\gamma_{1}}+\sqrt{\gamma_{2}}\right)^{1/2}}
{\sqrt[\leftroot{-3}\uproot{5}4]{\gamma_{1}\gamma_{2}}}\\
&\geq\frac{\varepsilon\sqrt[\leftroot{-2}\uproot{4}4]{\varsigma_{D}}}{4}\|\matr{\Lambda}(\matr{X}_{k-1})\||x_k^{*}-1|\sqrt{1+\frac{1}{{x_k^{*}}^2}}
\geq\frac{\varepsilon\sqrt[\leftroot{-2}\uproot{4}4]{\varsigma_{D}}}{4}\|\matr{\Lambda}(\matr{X}_{k-1})\|\|\matr{\Utwo}^{*}_{k}-\matr{I}_{2}\|. 
\end{align*}
\end{itemize}
Now we set $\iota_{D} = \min(\frac{3(1-4\varsigma_{D})\varepsilon}{4\sqrt{5}},\frac{\varepsilon\sqrt[\leftroot{-2}\uproot{4}4]{\varsigma_{D}}}{4})$. 
The proof is complete.
\end{proof}

\begin{proof}[Proof of \cref{lemma:global_rho_c-2-0}]
We prove that the inequality \cref{equation-condition-weak_c-3} is satisfied in two cases.
\begin{itemize}
\item If $h_k = h^{(U)}_{k}$, by \Cref{lemm:elemen_express_c}(i), we see that
\begin{align*}
g(\matr{X}_{k-1})-g(\matr{X}_{k})
&=h^{(U)}_{k}(0,0)-h^{(U)}_{k}(x_k^*,y_k^*) = \frac{1}{\alpha_{1}}(\alpha_{2}^2+\alpha_{3}^2)\\
&=\frac{1}{2}\|\partial h^{(U)}_{k}(0,0)\|\|(x_k^{*},y_k^{*})\|
\geq\frac{\varepsilon}{2}\|\matr{\Lambda}(\matr{X}_{k-1})\|\|\matr{\Utwo}^{*}_{k}-\matr{I}_{2}\|. 
\end{align*}
\item If $h_k = h^{(L)}_{k}$, by \Cref{lemm:elemen_express_c}(ii), we see that
\begin{align*}
g(\matr{X}_{k-1})-g(\matr{X}_{k})
&=h^{(L)}_{k}(0,0)-h^{(L)}_{k}(x_k^*,y_k^*) = \frac{1}{\beta_{1}}(\beta_{2}^2+\beta_{3}^2)\\
&=\frac{1}{2}\|\partial h^{(L)}_{k}(0,0)\|\|(x_k^{*},y_k^{*})\|
\geq\frac{\varepsilon}{2}\|\matr{\Lambda}(\matr{X}_{k-1})\|\|\matr{\Utwo}^{*}_{k}-\matr{I}_{2}\|. 
\end{align*}
\end{itemize}
Now we set $\iota_{LU}=\frac{\varepsilon}{2}$. The proof is complete.
\end{proof}

\begin{proof}[Proof of \cref{lemma:weak_rho_inequality_c-1}]
If $h_k=h^{(U)}_{k}$, we have
\begin{align*}
\|\matr{\Utwo}^{*}_{k}-\matr{I}_{2}\|^2=\frac{\alpha_2^2+\alpha_3^2}{\alpha_1^2}\geq \frac{\|\partial h^{(U)}_{k}(0,0)\|^2}{4\mathrm{M}_{\alpha}^2}\geq \frac{\varepsilon^2}{4\mathrm{M}_{\alpha}^2}\|\matr{\Lambda}(\matr{X}_{k-1})\|^2,
\end{align*}
where $\mathrm{M}_{\alpha}$ is a fixed positive constant always satisfying $|\alpha_1|\leq\mathrm{M}_{\alpha}$. 
The case $h_k=h^{(L)}_{k}$ is similar.
Now we prove the $h_k=h^{(D)}_{k}$ case. If $\varpi\in[0,\varsigma_{D})$, we have that 
\begin{align*}
\|\matr{\Utwo}^{*}_{k}-\matr{I}_{2}\|^2=\frac{5}{4}\geq\frac{5}{4}\frac{1}{4(\gamma_1^2+\gamma_2^2)}
\|\partial h^{(D)}_{k}(1,0)\|^2\geq\frac{5}{4}\frac{\varepsilon^2}{4\mathrm{M}_{0}^2}
\|\matr{\Lambda}(\matr{X}_{k-1})\|^2.
\end{align*}
The case $\varpi\in(\frac{1}{\varsigma_{D}},+\infty]$ is similar.
If $\varpi\in [\varsigma_{D},\frac{1}{\varsigma_{D}}]$, we have 
\begin{align*}
\|\matr{\Utwo}^{*}_{k}-\matr{I}_{2}\|^2\geq (1-x_k^{*})^2\geq \frac{\gamma_1^2(1-\varpi)^2}{\mathrm{M}_{0}^2\mathrm{M}_\varpi^2}=\frac{\|\partial h^{(D)}_{k}(1,0)\|^2}{4\mathrm{M}_{0}^2\mathrm{M}_\varpi^2}\geq\frac{\varepsilon^2}{4\mathrm{M}_{0}^2\mathrm{M}_\varpi^2}\|\matr{\Lambda}(\matr{X}_{k-1})\|^2,
\end{align*}
where $\mathrm{M}_\varpi$ is a fixed positive constant always satisfying $(1+\sqrt{\varpi})(1+\sqrt[\leftroot{-2}\uproot{4}4]{\varpi})\leq\mathrm{M}_\varpi$. 
We only need to set $\kappa^2$ to be the minimum of all the above corresponding positive constants. 
The proof is complete. 
\end{proof}

\section*{Acknowledgment} 
The authors would like to thank the three anonymous reviewers and the editor for their helpful suggestions and comments, which significantly improved the presentation of the article.

\bibliographystyle{siamplain}
\bibliography{Tensor_All}

\end{document}

%% file: tmp_BCD-G-SIMAX_header.tex
\title{Gradient based block coordinate descent algorithms for joint approximate diagonalization of matrices\thanks{Nov. 5th, 2021.
\funding{This work was supported in part by the National Natural Science Foundation of China (No. 11601371) and the GuangDong Basic and Applied Basic Research Foundation (No. 2021A1515010232).}}}

\author{Jianze Li\thanks{Shenzhen Research Institute of Big Data, The Chinese University of Hong Kong, Shenzhen, China (\email{lijianze@gmail.com}).}
\and Konstantin Usevich\thanks{Universit\'{e} de Lorraine, CNRS, CRAN, Nancy, France (\email{konstantin.usevich@univ-lorraine.fr}).} \and Pierre Comon\thanks{Univ. Grenoble Alpes, CNRS, Grenoble INP, GIPSA-Lab, France (\email{pierre.comon@gipsa-lab.fr}).}}

% Custom SIAM macro to insert headers
\headers{Gradient based block coordinate descent algorithms}{Jianze Li, Konstantin Usevich and Pierre Comon}

%% file: tmp_BCD-G-SIMAX_abstract.tex
\begin{abstract}
In this paper, we propose a gradient based block coordinate descent (BCD-G) framework to solve the joint approximate diagonalization of matrices defined on the product of the complex Stiefel manifold and the special linear group.
Instead of the cyclic fashion,
we choose the block for optimization in a way based on the Riemannian gradient.
To update the first block variable in the complex Stiefel manifold, we use the well-known line search descent method.
To update the second block variable in the special linear group, based on four different kinds of elementary transformations, we construct three classes: GLU, GQU and GU, and then get three BCD-G algorithms: BCD-GLU, BCD-GQU and BCD-GU.
We establish the global convergence and weak convergence of these three algorithms using the \L{}ojasiewicz gradient inequality under the assumption that the iterates are bounded. We also propose a gradient based Jacobi-type framework to solve the joint approximate diagonalization of matrices defined on the special linear group. Similar as in the BCD-G case, using the GLU and GQU classes of elementary transformations, we focus on the Jacobi-GLU and Jacobi-GQU algorithms, and establish their global convergence and weak convergence as well.
All the algorithms and convergence results in this paper also apply to the real case.
\end{abstract}

\begin{keywords}
blind source separation, joint approximate diagonalization of matrices, block coordinate descent, Jacobi-G algorithm, convergence analysis, manifold optimization
\end{keywords}

\begin{AMS}
49M30, 65F99, 90C30, 15A23
\end{AMS}